%% file: main.tex
\pdfoutput=1
\documentclass{article}

\usepackage{graphicx}
\usepackage[dvipsnames]{xcolor}

\usepackage{mathtools}
\usepackage{amsmath,amsthm,amssymb}
\usepackage{esint}
\usepackage{bbm}

\usepackage[
  margin=2.5cm,
  includefoot,
  footskip=30pt,
]{geometry}
\usepackage{layout}

\usepackage{enumitem}
\usepackage{hyperref}
\newcommand{\notoc}[1]{\texorpdfstring{#1}{}}
\hypersetup{hidelinks} 
\usepackage{cleveref}

\usepackage{comment}

\counterwithin*{equation}{section}

\newtheorem{theorem}{Theorem}[section] 
\newtheorem{lemma}[theorem]{Lemma}
\newtheorem{corollary}[theorem]{Corollary}
\newtheorem{proposition}[theorem]{Proposition}
\newtheorem{remark}[theorem]{Remark}

\newtheorem*{theorem*}{Theorem}
\newtheorem*{proposition*}{Proposition}

\theoremstyle{definition}
\newtheorem{definition}[theorem]{Definition}

\newcommand{\grad}{{\nabla}}

\newcommand{\R}{{\mathbb{R}}}
\newcommand{\E}{{\mathbb{E}}}

\DeclareMathOperator*{\argmin}{arg\,min}

\newcommand{\dd}{\mathop{}\!\mathrm{d}}
\newcommand{\1}{\mathbbm{1}}

\newcommand{\Barrel}[3][d]{\operatorname{F}_{#1} (#2,#3)}
\newcommand{\Slice}[3][d]{\operatorname{S}_{#1} (#2,#3)}

\title{Convex sets can have interior hot spots}
\author{Jaume de Dios Pont\footnote{ETHZ Department of Mathematics, \href{mailto:jaume.dediospont@math.ethz.ch}{jaume.dediospont@math.ethz.ch}}}
\date{}
\begin{document}
\maketitle

\abstract{
    The hot spots conjecture asserts that for any convex bounded domain $\Omega$ in $\mathbb R^d$, the first non-trivial Neumann eigenfunction of the Laplace operator in $\Omega$ attains its maximum at the boundary.
    We construct counterexamples to the conjecture for all sufficiently large values of $d$. 
    The construction is based on an extension of the conjecture from convex sets to log-concave measures.}

\setlength{\parskip}{.5em}

\section{Introduction}

\input{Parts/introduction}

\section{Proof of \notoc{\Cref{thm:counterexample}}}%
\label{sec:proof}
    \input{Parts/statement}

\section{Proof of \notoc{\Cref{thm:cylindrical_convergence}}}%
\label{sec:proof_cylinder}

    Through this section we will fix the convex pair \((\Omega,V)\), and denote by \(\Omega_d\) the set \(\Barrel{\Omega}{V}\). All the constants in this section may depend on \(\Omega\) and \(V\), but not on \(d\).
    
    \subsection{Basic Barrel properties}
        \input{Parts/basic_cylinder_properties}

    \subsection{\notoc{Behavior of \(\psi_d\) near \(t=0\)}}
        \input{Parts/boundary_l2}

\subsection{The limiting PDE on the interior}
    \input{Parts/heat_flow}

\section{Perturbations on a rectangle}%
\label{sec:perturbation}
    \input{Parts/perturbation}

\section{Wings}%
\label{sec:wings}
    \input{Parts/wings}

\bibliographystyle{alpha}
\bibliography{bibliography}

\crefalias{section}{appendix}
\renewcommand{\thetheorem}{A\arabic{theorem}} 
\renewcommand{\theproposition}{A\arabic{proposition}} 
\renewcommand{\thelemma}{A\arabic{lemma}} 
\renewcommand{\thesection}{A}
\renewcommand{\thesubsection}{A.\arabic{subsection}}
\setcounter{section}{0} 
\setcounter{subsection}{0}

\section{Step-by-step computations, Folklore, and Well-known results}
    \input{Parts/appendix}

\end{document}

%% file: Parts/introduction.tex
Let \(\Omega \subseteq \R^d\) be a bounded, connected domain, and let \(\phi_\Omega\) denote the first nontrivial Laplace eigenfunction in \(\Omega\) with Neumann boundary conditions. 
The \emph{Hot Spots} conjecture, posed by Ranch in 1974, asserts that for particularly \emph{simple}  sets \(\Omega\) the function \(\phi_\Omega\) attains its extreme values in \(\partial \Omega\). 
The conjecture is expected to hold for simply connected domains in \(\R^2\), and was expected to hold for convex domains in \(\R^d\) (see e.g~\cite{banuelos1999hot, Burdzy2005TheHS}, or~\cite{nickl2024consistent}). 
Although the natural classes of domains for which the conjecture should hold remain unclear, the conjecture is known to be false for the class planar domains with at least one hole~\cite{burdzy1999counterexample, Burdzy2005TheHS}.  

In \(\mathbb{R}^2\), the conjecture has been proven for several classes of domains, including all triangles~\cite{judge2020euclidean,judge2022erratum}; domains that are long and thin in a certain sense~\cite{banuelos1999hot,atar2004neumann}; and convex sets with one axis of symmetry~\cite{jerison2000hot,pascu2002scaling}.

In higher dimensions, the conjecture has been proven for several classes of domains. These include cylinder sets~\cite{kawohl1985rearrangements}, thin sets that are rotationally symmetric in all but one dimension~\cite{chen2019monotone}, and domains that naturally generalize the two-dimensional class introduced by Atar and Burdzy~\cite{atar2004neumann}, such as~\cite{yang2011hot, KennedyRohleder2024}. The strongest result without additional assumptions in higher dimensions is related to the \emph{Hot Spots Constant} \(S_d\), introduced by Steinerberger~\cite{steinerberger2023upper} as
\begin{equation}
	S_d:=\sup_{\substack{\Omega \subseteq \R^d\\ \Omega \text{ connected} \\ \partial \Omega \text{ smooth} }} \frac{\max_{x\in \Omega} \phi_\Omega(x)}{\max_{x\in \partial \Omega} \phi_\Omega(x)}.
\end{equation}
Steinerberger showed that the constant is uniformly bounded in all dimensions. The best known upper bounds, shown in~\cite{mariano2023improved}, include that \(S_d\le 5.0143\) in all dimensions, and that \(\limsup_{d\to \infty} S_d \le \sqrt{e}\).

This work focuses on convex sets in high dimensions. This version of the Hot Spots conjecture was first proposed by Kawohl~\cite{kawohl1985rearrangements}, with further refinements by Bañuelos and Burdzy~\cite{banuelos1999hot} and Jerison~\cite{jerison2019two}. The main result of this work is as follows:

\begin{definition}
	Let \(\Omega\subset \R^n\) be a convex domain, that is, a closed, bounded, convex subset of \(\R^n\) with nonempty interior. We denote by
	\(
		\lambda_\Omega
	\)
	the minimum nonzero Neumann eigenvalue of \(-\Delta\) in \(\Omega\). Let
	\(
		\phi_{\Omega}
	\)
	be any of the real Neumann eigenfunctions of \(-\Delta\)  in \(\Omega\) with eigenvalue \(\lambda_\Omega\), normalized with \(\frac{1}{|\Omega|}\int_{\Omega} |\phi_\Omega|^2 = 1 \). We say that \(\Omega\) has a spectral gap if \(\phi_\Omega\) is unique up to the sign.
\end{definition}

\begin{theorem}
	[The Hot Spots conjecture is false in high dimensions]
	\label{thm:counterexample}
	There exists a sequence \((\Omega_d)_{d\ge 1}\) of smooth, centrally symmetric domains with a spectral gap, with \(\Omega_d \subseteq \R^d\),  for which
	\begin{equation}
	\label{eq:counterexample_goal}
		\lim_{d\to \infty}  \frac{\max_{x\in \Omega_d} \phi_{\Omega_d}(x)}{\max_{x\in \partial \Omega_d} \phi_{\Omega_d}(x)} = \lim_{d\to \infty} \frac{\min_{x\in \Omega_d} \phi_{\Omega_d}(x)}{\min_{x\in \partial \Omega_d} \phi_{\Omega_d}(x)}  > 1.
	\end{equation}
	In particular, for \(d\) large enough, the maximum of \(\phi_{\Omega_d}\) is only attained on the interior.
\end{theorem}

This disproves the Hot spots conjecture for convex sets in sufficiently large dimensions, including the natural variations~\cite[Conjectures (HS1)-(HS3)]{banuelos1999hot}, and~\cite[Conjecture 2.5]{jerison2019two}. 
In comparison with existing positive results in the literature, the sets \(\Omega_d\) are most similar to the class of domains for which the conjecture was proven in~\cite{chen2019monotone}. The domains in ~\cite{chen2019monotone} are rotationally symmetric in all but one dimension. On the other hand, the sets \(\Omega_d\) in the counterexample are rotationally symmetric in all but two dimensions. Moreover, the sets \(\Omega_d\) in the counterexample are symmetric with respect to the operation \((x_1, \dots, x_i, \dots, x_d)\mapsto(x_1, \dots, -x_i, \dots, x_d)\), showing that a natural generalization of the hot spots conjecture for planar convex sets with two lines of symmetry cannot hold either.

The value of \(d_0\) after which the quotient is larger than one is not tracked through the proof, but one should be able to do so after straightforward modifications. The main required changes are quantifying the qualitative \emph{spectral gap} of certain constructions (e.g \Cref{prop:wings}) and quantifying the Fatou-type steps (such as at the end of the proof of (C.3)(a) of \Cref{thm:cylindrical_convergence}). However, a direct tracking of the argument would likely yield a dimension that is too large to be of any practical interest.

\subsection{Main ideas of the proof}

\textbf{Step 0: Posing a log-concave analog.} 
Dimension-free statements in convex analysis can be systematically extended from convex sets to log-concave measures. 
Let \(\Omega \subseteq \R^k\) be a convex set, and \(V\) a convex function in \(\Omega\).  
Consider sets of the form
\begin{equation}
	\Barrel{\Omega}{V}:={(x,w) \in \R^k\times \R^{d+1}, x\in \Omega, |w|\le C_d(1+ d^{-1/2}V(x))},
\end{equation} 
with \(C_d\) chosen appropriately (in our case \(C_d:= \frac{\sqrt{d}}2\), see \Cref{def:barrel}). 
In a certain sense, these sets approximate the log-concave measure \(\1_{\Omega}\exp(-V(x))\).  
For example, applying the Brunn-Minkowski inequality to sets of this form and sending \(d\to \infty\) yields the Prékopa-Leindler inequality. 
If the Hot Spots conjecture were true in all dimensions, one would expect a similar result to hold for log-concave measures. 
The proof begins by constructing such a log-concave analog. 

The first non-constant Neumann Laplace eigenfunction in \(\Omega\) is the minimizer of the Rayleigh quotient
\begin{equation}
	\frac{\int_{\Omega} |\nabla f(x)|^2 \1_{\Omega}dx}{\int_{\Omega} |f(x)|^2 \1_{\Omega}dx}.
\end{equation}
One can generalize the weight \(\1_{\Omega} dx\) to a general log-concave measure \(\exp(-V) dx\). For general \(V\), the critical points of the Rayleigh quotient
\begin{equation}
	\frac{\int_{\mathbb R^k} |\nabla f(x)|^2 e^{-V(x)}dx}{\int_{\mathbb R^k} |f(x)|^2 e^{-V(x)}dx}
\end{equation}
are the eigenfunctions \(\phi_{\Omega, V}\) of  \(L_{\Omega,V}:=-\Delta + \nabla V\cdot \nabla\) with Neumann boundary data in \(\partial \Omega\). 
The {log-concave} extension of the hot spots conjecture would assert that the first non-constant eigenfunction of \(L_{\Omega, V}\) attains its maximum in \(\partial \Omega\). 
This extension is false\footnote{Together with A. Salim, I had previously attempted to prove this log-concave extension of the conjecture, to solve the hot spots conjecture in the positive. The strategy was to study the operator \(-\epsilon \Delta+\nabla V \cdot \nabla\) using a homotopy from \(\epsilon = 0\) to \(\epsilon = 1\). The reason this approach fails (beyond the fact that the conjecture is false) is that the flow lines of the vector field \(\nabla V (x)\) can be quite complex even when \(V\) is convex.
A form of this idea can be seen in \Cref{prop:wings}.}. This suggests setting \(\Omega_{d}:=\Barrel{\Omega}{V}\) for some well-chosen \(\Omega, V\).\\

\noindent \textbf{Step 1: Understanding the log-concave analog.} (\Cref{thm:cylindrical_convergence}) The proof then proceeds by studying \(\phi_{\Omega_d}(x,w)\). First, one shows that the first eigenvalue of \(\Omega_d\) converges to the first eigenvalue of \(L_{\Omega, V}\).  Then one shows that the functions \(\phi_{\Omega_d}(x,w)\) are radial in the \(w\) variable. This allows us to set the ansatz \(\phi_{\Omega_d}(x,w) = \psi_d(x, 1-\frac{2 |w|^2}{d})\). This reduces the number of variables from \(\text{dim}(\Omega)+d+1\) to \(\text{dim}(\Omega)+1\). 

Most of the volume of \(\Barrel{\Omega}{V}\) is found in a \(\sim d^{-1/2}\) neighborhood of the set \(|w| = \frac{\sqrt{d}}{2}\). If the first eigenfunction of \(L_{\Omega,V}\) is unique, one expects the functions \(\psi_{d}(x,r)\) to converge in \emph{some} sense to \(\phi_{\Omega, V}(x)\). The change of variables \((x,w)\mapsto (x,t) = \left (x, 1-\frac{2 |w|^2}{d}\right)\) sends most of the volume of  \(\Barrel{\Omega}{V}\)  to a \(d^{-1}\)-neighborhood of \(t=0\). In particular, one expects \(L^2\) estimates on \(\phi_{\Omega,V}\) to give information about \(\psi_d\) only near \(t=0\).  This is the case: \(\psi_d(x,0)\) converges in \(L^2\) to the function \(\psi_{\Omega,V}\) as \(d\to \infty\). The convergence is upgraded to uniform convergence using the uniform Harnack inequalities of Li and Yau~\cite{li1986parabolic,wang2006dimension}.  

The nonlinear change \((x,t) = \left (x, 1-\frac{2 |w|^2}{d}\right)\) is useful when \(t>0\). Changing variables in the equation \(-\Delta \phi_{\Omega_d} = \lambda_{\phi_{\Omega_d}}\) one formally obtains (see eq.~\eqref{eq:psi_com_pde}) the equation
\begin{equation}
	\partial_t \psi_d = \frac 1 8 (\Delta + \lambda_d) \psi_d +O(d^{-1}),
\end{equation}
where the term \(O(d^{-1})\) contains first and second derivatives of \(\psi_d\). This will imply that functions \(\psi_d\) converge uniformly to a function \(h\) in \(\Omega\times [0,1]\), that solves
\begin{equation}
\label{eq:parabolic_intro}
	\begin{cases}
		\partial_t h = \frac 1 8 (\Delta + \lambda_{\Omega, V}) h \\
		h(x,0)= \psi_{\Omega,V}\\
		\partial_{\vec n} h(x,t)=0 \text{ in } \partial \Omega \times [0,1]
	\end{cases}
\end{equation}
The change of variables from \(\phi_{\Omega_d}\) to \(\psi_d\) maps (as \(d\to \infty\)) the boundary \(\partial \Omega_d\) to the set \((\Omega \times \{0\})\cup \partial \Omega \times [0,1]\). This is precisely the parabolic boundary of the PDE~\eqref{eq:parabolic_intro}. The goal now is reduced to finding such \(h\) that does not satisfy the parabolic maximum principle. When \(\kappa >0\), solutions to \(\partial_t h = \frac 1 8 (\Delta + \kappa) h\) do not satisfy the maximum principle for generic initial data. The challenge lies in finding initial data of the form \(\phi_{\Omega, V}\) for some convex potential \(V\).

Equation~\eqref{eq:parabolic_intro} is a parabolic equation obtained as a limit of elliptic equations. This is a well-documented phenomenon. It is, for example, a key step of the proof of the geometrization conjecture~\cite[Section 6]{perelman2003entropy}. See for example~\cite{davey2018parabolic} for a more recent approach to similar questions and historical remarks.

\noindent \textbf{Step 2: Preparing suitable initial data to \(h\).} (\Cref{prop:perturbation}) In order to find suitable initial data to ~\eqref{eq:parabolic_intro} we resort to a two-step perturbative argument. As a first step, we consider a small potential \(\epsilon \cdot V_0\) on a rectangle \(\Omega = R := [-\pi/2,\pi/2]\times [-1,1]\). Using first-order perturbation methods, this potential is seen to have eigenfunction \(\sin(x) + \epsilon \beta(x,y) +O(\epsilon^2)\) for some \(\beta(x,y)\). An explicit construction shows that one can prescribe the value of \(\beta(\pm \pi/2, y)\) while keeping \(V\) convex. At this stage, one still has no control over the maximum of \(h\), but one expects it to be at the boundary point \((x,y),t = (\pi/2, 0), 1\).

\noindent \textbf{Step 3: Shielding the maximum.} (\Cref{prop:wings}) The remaining step is to extend \(R\) to \(R_m := [-\pi/2-m, \pi/2+m]\times[-1,1]\), in such a way that the maximum of \(h\) stays near \((\pi/2, 0), 1\). By extending \(V\) to a \(V_m\) whose gradient is very large in \(R_m \setminus R\), one can  approximate the equation \((-\Delta + \nabla V \nabla - \lambda_{\Omega, V} )\phi_{R_m,V+,} = 0\) by \(\nabla V \nabla\phi_{\Omega,V} = 0\). This is \emph{transporting} the value of \(\phi_R\) in \(R_m\setminus R\) by the flow of \(\nabla V_m\). Choosing a suitable \(V_m\) ensures that \(\phi_R(m+\pi/2, y) \le \phi_R(\pi/2, y)\), where equality does not hold for all \(y\). By a suitable strong maximum principle argument, for \(t>0\) the maximum cannot be at \(x= m+\pi/2\) and has to be on the interior.

\subsection*{Notation} Let \((\Omega_n)_{n\ge 0}\) be a sequence of sets in \(\R^d\), with continuous functions \(f_n:\Omega_n \to \R\).

The expression \emph{the functions \(f_n\) are equicontinuous} will denote that for any \(\epsilon>0\) there is an \(\delta>0\) such that for any \(n\) and \(x,y\in \Omega_n\) such that \(|x-y|<\delta\), \(f_n(x)-f_n(y)<\epsilon\) (note that their domains need not be the same). The expression \emph{the functions \(f_n\) in \(\Omega_n\) converge uniformly to \(f\) in \(\Omega\)} will be used to denote that the sets \(\Omega_n\) converge to \(\Omega\) in the Hausdorff sense and  for any continuous function \(g\) in \(\R^d\) such that \(g|_{\Omega}=f\), one has \(\lim_{n\to \infty} \|f_n-g\|_{L^\infty(\Omega_n)} \to 0\).

The expression \emph{the functions \(f_n\) converge to \(f\)} up to a sign means that there exist signs \(\epsilon_n=\pm1\) such that \(\epsilon_n f_n\) converges to \(f\). \\ 

\textbf{Acknowledgements:}  I am grateful to Marco Badran, Enric Florit-Simon,  Svitlana Mayboroda, and Adil Salim for many enlightening conversations through this work.  This work was started after visits to S. Steinerberger at the University of Washington, while I was employed at Microsoft Research under the mentorship of Adil Salim and Jerry Li. I am extremely grateful for their guidance and encouragement during that time. 

The work was partially developed during the Hausdorff Institute Dual Trimester Program  \emph{
``Synergies between modern probability, geometric analysis and stochastic geometry''} funded by the Deutsche Forschungsgemeinschaft (DFG, German Research Foundation) under Germany's Excellence Strategy – EXC-2047/1 – 390685813. I thank the Hausdorff Institute, as well as the organizers of the program, for their hospitality.

%% file: Parts/statement.tex
This section aims to prove \Cref{thm:counterexample}, deferring the proofs of the more technical propositions (\Cref{thm:cylindrical_convergence}, \Cref{prop:perturbation}, and \Cref{prop:wings}) to the subsequent sections. These three results are motivated at an informal level in the current section, and only fully developed in Sections~\ref{sec:proof_cylinder}-\ref{sec:wings}. 
\subsection{From convex sets to log-concave measures}
The key step to building counterexamples witnessing \Cref{thm:counterexample} is to lift them from an extension of the Laplace equation to a Langevin equation with a potential:

\begin{definition}
	A convex pair in \(\R^n\) is a pair \((\Omega, V)\), where \(\Omega\) is a bounded convex domain in \(\R^n\), and \(V: \Omega \to \R\sqcup \{+\infty\}\) is a convex function. Define \(Z_{\Omega,V}:=\int_{\Omega} \exp(-V(x)) dx\), and the measure \(\mu_{\Omega,V}:= Z_{\Omega,V}^{-1}\exp(-V(x)) dx\). If \(V\) is zero it will be omitted, leaving \(\mu_{\Omega}=\mu_{\Omega,0}\).

	We define \(\lambda_{\Omega,V}\) to be the minimum, over non-zero functions \(f\) with mean zero, of the Rayleigh quotient
	\begin{equation}
		R_V(f):=\frac{\int_\Omega |\grad f|^2 d\mu_{\Omega,V}(x)}{\int_\Omega |f|^2 d\mu_{\Omega,V}(x)}.
	\end{equation}
	By the Rayleigh quotient formula, \(\lambda_{\Omega,V}\) is the lowest eigenvalue of the operator
	\begin{equation}
		L_{\Omega,V}f:=- \Delta f +\grad V \cdot \grad f
	\end{equation} 
	with Neumann boundary data. Let \(\phi_{\Omega,V}\) denote the any of the eigenfunctions with eigenvalue \(\lambda_{\Omega,V}\) with the normalization \(\|\phi_{\Omega,V}\|_{L^2(\mu_{\Omega,V})} = 1 \). The pair \((\Omega, V)\) will be said to have a spectral gap if \(\phi_{\Omega, V}\) is unique up to sign. If \(\Omega\) has a smooth boundary, and \(V\) extends to a smooth function in a neighborhood of \(\Omega\), the pair \((\Omega, V)\) is a \emph{smooth} convex pair.
\end{definition}

\begin{remark}
	Even when \((\Omega, V)\) has a spectral gap, the function \(\phi_{\Omega, V}\) is still only defined up to a sign. The convention chosen through the paper (which is specified in each case) is in general that
	\begin{equation}
		\int_{\Omega} \phi_{\Omega,V}(x) \cdot x_1 d\mu_{\Omega,V}(x) >0.
	\end{equation}
\end{remark}

The definition of \(\phi_{\Omega, V}\) generalizes that of \(\phi_\Omega\), in the sense that \(\phi_{\Omega,0}=\phi_\Omega\). In the converse direction, one can \emph{simulate} \(\phi_{\Omega, V}\) by designing a convex set \(\Barrel{\Omega}{V}\) which we call a barrel domain.

\begin{definition}[Barrel domain]\label{def:barrel}
	Let \((\Omega, V)\) be a convex pair in \(\R^n\) such that \(\max_{x\in \overline \Omega}V(x)\le 0\), and let \(d\ge 2 \|V\|_{L^\infty(\Omega)}\). Let
	\begin{equation}
		\rho_d(V)(x):= \frac{1}{2}\left(\sqrt{d}-\frac{V(x)}{\sqrt{d}}\right)
	\end{equation}
 The barrel \(\Barrel{\Omega}{V} \subseteq \R^{d+n+1}\) is the convex domain
	\begin{equation}
  	\Barrel{\Omega}{V} := \left\{(x,w)\in \Omega\times\R^{d+1}, |w|\le \rho_d(V)(x)\right\}.
\end{equation} 
We denote \(S_d(\Omega,V)\) the \emph{slice} domain
\begin{equation}
	S_d(\Omega,V):=\{(x,\sqrt{d}/2-|w|) , (x,w)\in \Barrel{\Omega}{V} \}.
\end{equation}
To the slice domain we will associate the slice potential \(W^s_d(x,r):= -\sqrt{d}\log(1- 2/\sqrt{d}r)\).
\end{definition}

\begin{remark}
	 The potential \(W^s_d(x,r)\) is such that the probability measure \(\mu_{\Slice{\Omega}{ V},\sqrt{d}W^s_d}\) is the pushforward of the uniform probability \(\mu_{\Barrel{\Omega}{V},0}\) by the map \((x,w)\mapsto (x, \sqrt{d}/2-|w|)\).
\end{remark}

\begin{remark}
	The name of Barrel domain stems from the simplest example of such domain: when \(d=1\), and \(\Omega\) is an interval, the set \(\Barrel[1]{\Omega}{V}\) looks like a barrel. The slice domain is a \emph{slice} of this barrel, displaced so that its minimum is near zero:
	\begin{center}
	\raisebox{-0.5\height}{
	\includegraphics[width=.18\textwidth]{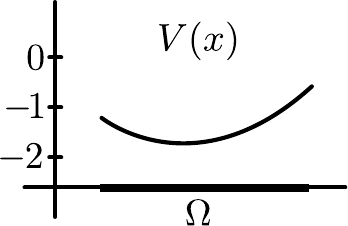}}
	\(\qquad \to\qquad \)
	\raisebox{-0.5\height}{
	\includegraphics[width=.15\textwidth]{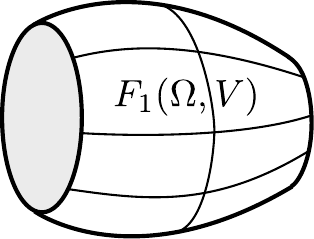}}
	\(\qquad \to\qquad \)
	\raisebox{-0.5\height}{
	\includegraphics[width=.2\textwidth]{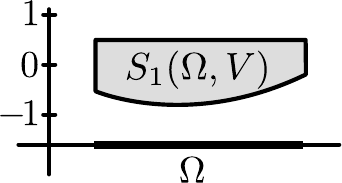}}
	\end{center}
	When \(\Omega\) is a rectangle (the case of interest) the set \(F_1(\Omega, V)\) is already \(4-\)dimensional. In the case of \Cref{thm:counterexample}, the set \(\Omega\) will instead be a 2-dimensional rectangle.
\end{remark}

The factor \(\frac{1}{2}\) in the definition of Barrel domains is arbitrary and chosen for convenience later in the proof.  The ground eigenfunctions \(\phi_{\Barrel{\Omega}{V}}\) are closely related to the ground eigenfunctions of \(\phi_{V,\Omega}\) by the following result:

\begin{theorem}[Eigenfunctions on barrel domains]
    \label{thm:cylindrical_convergence}
    \label{prop:cylindrical_convergence} 
    Let \((\Omega, V)\) be a smooth convex pair with \(V\ge 0\) a spectral gap, and \(\lambda_{\Omega,V}\le 2\). Let \(\Omega_d:=\Barrel{\Omega}{V}\), then:
    
 \begin{enumerate}[label ={(C.\arabic*)}]
     \item \label{item:eigenvalues_converge} The ground eigenvalues \(\lambda_{\Omega_d}\) converge to \(\lambda_{\Omega,V}\) as \(d\to \infty\), and for \(d\) large enough, the sets \(\Omega_d\) have a spectral gap. 
     \item \label{item:eigenfunctions_radial}  For all large enough \(d\) the ground eigenfunction  \(\phi_{\Omega_d}(x,w)\) of \((\Omega_d,0)\) is radial in the second coordinate, and thus can be written in the form  \(\phi_{\Omega_d}(x,w)=: \psi_d(x,1-(|w|/\rho_d(V)(x))^2)\). 
     \item \label{item:limit} The functions \(\psi_{d}(x,t)\) converge up to a sign weakly in \(L^2(\Omega \times[0,1])\) to a sequence \(\psi_\infty\) that satisfies:
     \begin{enumerate}
        \item \(\psi_\infty(x,0) = \phi_{\Omega, V}(x)\) 
         \item  \(\partial_t \psi = \frac{1}{8}(\Delta \psi + \lambda_0 \psi)\) on \(\Omega^\circ \times (0, \tau)\). 
         \item Neumann boundary conditions on \(\partial\Omega \times (0, 1)\) 
     \end{enumerate}
      \item \label{item:convergence} The functions \(\psi_d(x,t)\) converge up to a sign in the \(C^0\) topology to the function \(\psi_\infty\).
\end{enumerate}
\end{theorem}

The proof of \Cref{thm:cylindrical_convergence} is deferred to \Cref{sec:proof_cylinder}. The main ideas of the proof are as follows.
\begin{enumerate}
	\item By symmetrization there exists an eigenbasis such that the eigenfunctions in \(\Barrel{\Omega}{V}\) are either radial in the second variable or have mean zero on spheres in the second variable. The Poincaré constant in balls of radius at most \(\sqrt{d}\) in \(\mathbb{R}^{d+1}\) is less than \(1\), which rules out the second scenario.
	\item Most of the volume of \(\Barrel{\Omega}{V}\) is held in the shell consisting of \((x,w)\) for which \(||w|- \sqrt{d}/2|\lesssim \frac 1{\sqrt{d}}\). This shell only has thickness \(d^{-1/2}\), and one expects eigenfunctions to be roughly constant in the \(|w|\) variable at scale \(d^{-1/2}\). If the functions \(\phi_{\Omega_d}(x,w)\) depended only on the \(x\) variable, they would minimize the Rayleigh quotient
	\begin{equation}
		\frac{\int_{\Omega}|\grad_x \phi_{\Omega_d}(x)|^2 \left(\sqrt d - \frac{V(x)}{\sqrt d}\right)^d dx}{\int_{\Omega}|\phi_{\Omega_d}(x)|^2 \left(\sqrt d - \frac{V(x)}{\sqrt d}\right)^d dx}.
	\end{equation}
	Given that
	\begin{equation}
		 \sqrt{d}^{-d}\left(\sqrt d - \frac{V(x)}{\sqrt d}\right)^d \to \exp(-V(x)),
	\end{equation}
	one would expect \(\|\phi_{\Omega_d}-\phi_{\Omega,V}\|_{\mu_{\Barrel{\Omega}{V}}}\to 0\). Since most of the volume of \(\Barrel{\Omega}{V}\) is held in a thin shell at the boundary, one expects this convergence to hold in the \(L^\infty\) sense in that shell, but not necessarily on the interior.

	\item Let \(\psi^{com}_d(x,t)\) be defined implicitly by \(\phi_{\Omega_d}(x,w) = \psi^{com}_d\left (x,1-\left (\frac{|w|}{\sqrt{d}/2}\right)^2\right)\). The functions \(\psi_d, \psi^{com}_d\) are related by continuous changes of variables,
	\begin{equation}
	\label{eq:psi_com}
		\psi_d(x,t)=\psi_d^{com}(x, (1-d^{-1}V(x))^2t+2 d^{-1}V(x) - d^{-2}V(x)^2).
	\end{equation} These changes of variables converge to the identity as \(d\to \infty\).
	The functions \(\psi^{com}_d(x,t)\) (See \Cref{sub:step_by_step} for a step-by-step computation of this change of variables) satisfies the PDE
	\begin{equation}
	\begin{aligned}
	\label{eq:psi_com_pde}
		-\Delta_x \psi^{com}_d  +  
		   8 \partial_t \psi^{com}_d = \lambda_{\Omega_d}  \psi^{com}_d
		  -
		   \underbrace{\frac{16}{d} \left (\frac{|w|}{\sqrt{d}/2}\right)^2 \Delta_t \psi^{com}_d + 
		\frac{8}{d} \partial_t  \psi^{com}_d}_{\text{Formally }\to 0}
	\end{aligned}
	\end{equation}
	with Neumann boundary conditions on \([0,1]\times\partial\Omega\). 

	\item By Li-Yau's uniform Harnack inequality, the eigenfunctions should be equicontinuous, and thus have convergent subsequences in \(C^0\) satisfying all of the properties.
\end{enumerate}

\Cref{thm:cylindrical_convergence} motivates the definition of the Heat extension function:

\begin{definition}
	Let \((\Omega, V)\) be a convex pair in \(\R^n\) with a spectral gap. Let \(h_{\Omega_V}:\Omega \times [0,1]\) be the unique continuous function (up to a sign) satisfying 
	\begin{equation}
	\label{eq:h_extension}
		\begin{cases}
		    	h_{\Omega,V}(x,0) = 
		    		\phi_{\Omega, V}(x)
		    \\
         		\partial_t h_{\Omega,V} = 
         			\frac{1}{8}(\Delta_x h_{\Omega,V}+ \lambda_{\Omega, V} h_{\Omega,V})
         		& 
         		\text{ on } \Omega^\circ \times (0, 1)
         	\\
         		\partial_{\vec n} h_{\Omega,V} = 0 
         		& 
         		\text{ on } \partial \Omega \times (0,1).
		\end{cases}
	\end{equation}
	The \emph{parabolic boundary} of \(\Omega\times [0,1]\), denoted by \(\partial_h \Omega\), is the set \(\partial \Omega \times [0,1] \cup \Omega \times \{0\}\).
\end{definition}

Once \Cref{thm:cylindrical_convergence} is shown, \Cref{thm:counterexample} reduces to showing
\begin{proposition}[No parabolic maximum principle]
\label{prop:heat_flow_spot}
	There exists a convex pair \((\Omega,V)\) with \(\|V\|_\infty <\infty\) and a spectral gap such that \(h_{\Omega,V}\) does not attain its maximum at \(\partial_h \Omega\). Moreover, the pair can be chosen so that \(\Omega\) is a rectangle, \(V\) is a symmetric function, and \(\phi_{\Omega, V}\) is antisymmetric in \(y\) and symmetric in \(x\).
\end{proposition}

The forcing term \(\lambda_{\Omega_d} h_{\Omega,V}\) to the heat equation~\eqref{eq:h_extension} has a positive sign. This arises from eq.~\eqref{eq:psi_com_pde}), and is central to the counterexample: The equation \(\partial_t f = \Delta f + \kappa f\) with Neumann boundary data satisfies the parabolic maximum principle if and only if \(\kappa\le 0\). \Cref{prop:heat_flow_spot} states that there are initial data of the form \(\phi_{\Omega, V}\) for which this forced heat flow does not satisfy the parabolic maximum principle for small times.

\begin{proof}
	[Proof that \Cref{prop:heat_flow_spot} and \Cref{thm:cylindrical_convergence} together imply \Cref{thm:counterexample}]

	Let \((\Omega, V)\) be a convex pair with a spectral gap such that the heat extension takes the maximum on the interior, that is such that for some \(\epsilon>0\),
	\begin{equation}
		\frac{\|h_{\Omega,V}\|_{L^\infty(\Omega\times [0,1])}}{\|h_{\Omega,V}\|_{L^\infty(\partial_h \Omega)}} \ge 1+4\epsilon
	\end{equation}
	One can approximate \((\Omega, V)\) with smooth counterparts, (see e.g. \Cref{lem:fokker_plank_stable} for a self-contained argument regarding approximability): There is a smooth convex pair \((\tilde \Omega, \tilde V)\) with a spectral gap such that 
	\begin{equation}
		\frac{\|h_{\tilde\Omega,\tilde V}\|_{L^\infty(\tilde\Omega\times [0,1])}}{\|h_{\tilde\Omega,\tilde V}\|_{L^\infty(\partial_h \tilde \Omega)}} \ge 1+3\epsilon.
	\end{equation}

	Let \(\phi_d:= \phi_{\Barrel{\tilde \Omega}{\tilde V}}\).
	The map \((x,w)\mapsto(x,1-(|w|/\rho_d(V)(x)^2)\) is a surjective map from \(\Omega_d\) to \(\Omega \times [0,1]\) that maps \(\partial \Omega_d\) to \(\partial_h \Omega\). Applying \Cref{thm:cylindrical_convergence}, (C.4)
	\begin{equation}
		\lim_{d\to \infty}
		\frac{\|\phi_d\|_{L^\infty( \phi_{\Barrel{\tilde \Omega}{ \tilde V} })}}{\|\phi_d\|_{L^\infty(\partial \phi_{\Barrel{\tilde \Omega}{ \tilde V} })}}
		=
		\frac{\|h_{\tilde\Omega,\tilde V}\|_{L^\infty(\tilde\Omega\times [0,1])}}{\|h_{\tilde\Omega,\tilde V}\|_{L^\infty(\partial_h \tilde \Omega)}}
		\ge
		1+2 \epsilon.
	\end{equation}
	By approximating each \(\Barrel{\Omega}{V}\) with smooth sets, for all \(d\) large enough there exists a convex smooth domain \(T_d\) with a spectral gap such that
	\begin{equation}
		\frac
		{\|\phi_d\|_{L^\infty( \phi_{T_d})}}
		{\|\phi_d\|_{L^\infty(\partial \phi_{T_d})}} 
		\ge 1+\epsilon.
	\end{equation}
\end{proof}

Pictorially, the relationship between \(\phi_{\Barrel{\Omega}{V}}\) and \(h_{\Omega,V}\), and their domains, is as follows:

\begin{center}
	\raisebox{-0.5\height}{
	\includegraphics[width=.2\textwidth]{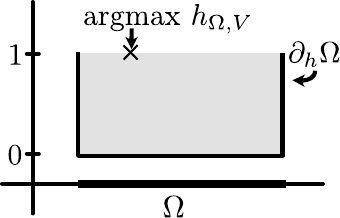}}
	\(\qquad \to\qquad \)
	\raisebox{-0.5\height}{
	\includegraphics[width=.21\textwidth]{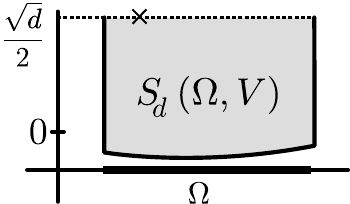}}
	\(\qquad \to\qquad \)
	\raisebox{-0.5\height}{
	\includegraphics[width=.18\textwidth]{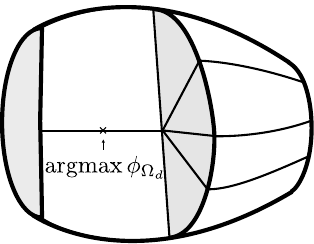}}
\end{center}

\begin{remark}
	The heat equation satisfies the parabolic maximum principle, which apparently could prevent \Cref{prop:heat_flow_spot} from holding. The positive term \(\lambda_0 h_{\Omega, V}\) on the right-hand side prevents the existence of such a maximum principle. 
\end{remark}

\subsection{Perturbation of a Laplace eigenfunction via a small potential}

The construction of \((\Omega,V)\) witnessing \Cref{prop:heat_flow_spot} will follow from two perturbation results around the convex set 
\begin{equation}
\label{eq:domain}
	R:=[-\pi/2,\pi/2] \times [-1,1],
\end{equation} which has unique first eigenfunction \(\phi_{R}(x,y)= \pm \sqrt{2/\pi}\sin(x)\) and eigenvalue \(\lambda_R=1\). We will take the convention \(\phi_{R}(x,y)= +\sqrt{2/\pi}\sin(x)\). The first perturbation result is

\begin{proposition}
\label{prop:perturbation}
	Let \(q(y):[-1, 1]\to \mathbb C\) be a smooth symmetric function with \(q'(\pm 1) = 0\). There exists a smooth function \(\beta:R\to \R\) such that
	\begin{enumerate}
		\item \(\beta(x,y)\) is symmetric in \(y\).
		\item \(\beta(x,y)\) is antisymmetric in \(x\).
		\item \(\beta(\pi/2, y) - q(y) = C_q\) is constant.
	\end{enumerate}

	and a smooth convex potential \(V_q:R \to \R_{\le 0}\) parametrized  such that as \(\epsilon \to 0\) the pair \((R,\epsilon V)\) has a spectral gap and
	\begin{equation}
	\label{eq:approx_limit}
		\begin{cases}
			\lambda_{R, \epsilon V_q} = 1 + O(\epsilon)\\
		\|\phi_{R, \epsilon V} - \sqrt{2}\sin(x)-\epsilon \beta \|_\infty = o(\epsilon)
		\end{cases}.
	\end{equation}
\end{proposition}
 
The proof of \Cref{prop:perturbation} is deferred to \Cref{sec:perturbation}. It follows from applying perturbation theory to the PDE \((-\Delta +\epsilon \nabla V \cdot \nabla )\phi_\epsilon = \lambda_{\epsilon}\phi_\epsilon\). If \(\beta = \frac{\partial \phi_\epsilon}{\partial \epsilon}\),  it must satisfy
\begin{equation}
\label{eq:pert_first}
	-\Delta \beta + \nabla V \cdot \nabla \phi_{R} = \lambda \beta + \mu \phi_R,
\end{equation}
for some potential \(V\). Given \(\beta\), one can construct a non-convex \(C^2\) potential \(V_0\) that ensures~\eqref{eq:pert_first} holds. Then one considers \(V= V_0 + M (x^2+y^2)\) for a large constant \(M\). The term \(Mx^2\) will give rise to the constant \(C_q\), and the term \(My^2\) has no effect on~\eqref{eq:pert_first}, since \(\nabla V \cdot \nabla \phi_{R} = \sqrt{2/\pi} \cos(x) \partial_x V \).

In the application of \Cref{prop:perturbation}, \(q(y)\) will be chosen to be any function with the properties below.
\begin{lemma}\label{def:hotspots_proposal}
 	There exists a smooth, symmetric function \(q(y):[-1,1]\to [-1,1]\) satisfying 
\begin{equation}
\label{eq:hotspots_proposal}
	\begin{cases}
		q(y) = 1 & \text{for } |y| \le \frac 13 \\
		q(y) = -1 & \text{for } \frac 23 \le y \le \frac 34\\
		q(1) = 0\\
		q'(y)\ge 0 & \text{for }  \frac  3 4 \le y \le 1 \\ q'(1)=0
	\end{cases}
\end{equation}
and such that \(H_q(y,t)\), the Neumann heat flow of \(q\) in \([-1,1]\) satisfies \(H_q(0,t)-H_q(1,t)>0\) for all  \(t\in [0,1]\).
 \end{lemma} 

One of these possible \(q(y)\), with an exaggerated sketch of the function \(\psi_{R, \epsilon V_q}\) is represented below.

\begin{center}
	\raisebox{-0.5\height}{
	\includegraphics[width=.2\textwidth]{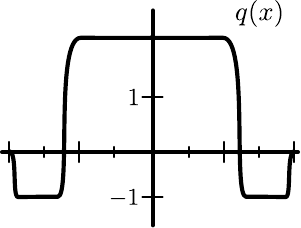}}
	\(\qquad \to\qquad \)
	\raisebox{-0.5\height}{
	\includegraphics[width=.35\textwidth]{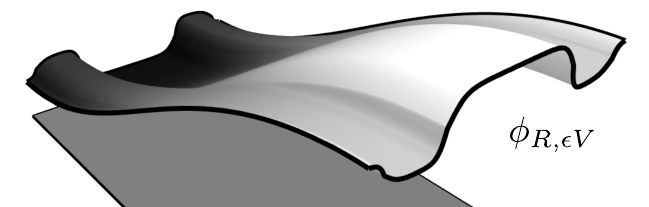}}
\end{center}

\begin{proof}[Proof of \Cref{def:hotspots_proposal}]
	Let \(b\) be a nonnegative symmetric decreasing bump function, such that \(b(x)=1\) for \(|x|\le \frac 13\), such that for \(x\in[0,1]\), \(b(x)= 1-b(1-x)\).  Then, for \(\delta\) small enough, the function
	\begin{equation}
		q_\delta(y) := 2 b(y)-1 + b(\delta^{-1}(|y|-1))
	\end{equation} satisfies the hypotheses.
\end{proof}

\subsection{\notoc{The operator \(-\Delta +  \grad V \cdot \grad\) as a perturbation of \(\grad V \cdot \grad\)}}

The pair \((R, \epsilon V_{q})\) (for \(\epsilon\) small enough)  still does not lead to a counterexample, one can show that the maximum of its heat extension is attained at or very close to \(x=\pi/2, y=0\). To produce an example one needs to extend \(R\) to \emph{shield} this counterexample at the boundary. We will do so by extending \(R\) to the set \(R_m:= [-\pi/2-m, \pi/2+m]\times [-1,1]\), for large values of \(m\). The potential \(V\) will be extended to a potential \(V_m\), such that in \(R_m\setminus R\), \(|\nabla V_m|\approx m^2\). In this case, the operator \(-\Delta +  \grad V_m \cdot \grad\) will be treated as a perturbation of \(\grad V_m \cdot \grad\).

\newcommand{\Wing}[3][m]{\operatorname{Wing}_{#1}({#2},{#3})}
\newcommand{\Trans}[2][g]{\operatorname{Ex}_{#1}(#2)}

\begin{definition}[Transport]\label{def:transport}
	 Let \(g(x,y)\) be a \(C^2([-1,1]\times[0,1])\) function such that:
	\begin{enumerate}
		\item For each \(x\), \(g(x,y)\) is convex in \(y\).
		\item \(g(x,y)\) is symmetric in the \(y\) variable.
		\item There is \(\delta>0\) such that \(g(x,y)=0\) for all \(x<\delta\).
	\end{enumerate}

	Let \(F(x_0,y_0): [0,1]\times[-1,1] \to [-1,1]\) be the value of \(\gamma(0)\) where \(\gamma\) is the solution to the ODE
	\begin{equation}
		\begin{cases}
			\gamma'(x) = \partial_y g(\gamma x))\\
			\gamma(x_0)=y_0
		\end{cases}
	\end{equation}

	Let \(f: R\to \mathbb R\) be a continuous function. Then we define 
	\begin{equation}
		\Trans f (x,y):= f(\pi/2, F(x,y)).
	\end{equation}

	In other words, the function \(\Trans f (x,y)\) takes the value of \(f(\pi/2, *)\), and \emph{extends} it by the vector field \((1, \partial_x g)\). In particular, if \(f\) is differentiable, \(\nabla \Trans  f \cdot (1, \partial_x g) = 0\).
\end{definition}

One can extend \((R,V)\) to a longer rectangle using the transport:

\begin{proposition}\label{prop:wings}
	Let \((R, V)\) be a convex pair, with \(R\) as in \eqref{eq:domain} such that \(\phi_{R, V}\) is antisymmetric in the first variable and \(V\) is Lipschitz. Let \(g\) satisfy the hypothesis of \Cref{def:transport}.

	Then there exists a sequence of convex pairs \((R_m, \Wing{V}{g})\), with \(R_m:=[-\pi/2 -m, \pi/2+m]\times [-1,1]\) such that, as \(m\to \infty\):
	\begin{enumerate}
		\item The ground eigenvalues \(\lambda_{R_m, \Wing{V}{g}}\) converge to \(\lambda_{R,V}\), and  for \(m\) large enough pair \((R_m, \Wing{V}{g})\) has a spectral gap.
		\item The functions \(\phi_{R_m, \Wing{V}{g}}(x,y)\) are antisymmetric in \(y\).
		\item In \(R\), the functions \(\phi_{R_m, \Wing{V}{g}}\) converge uniformly to \(\phi_{R,V}\)
		\item The functions \(T_m(x,y):=\phi_{R_m, \Wing{V}{g}}(m x+\pi/2, y)\) converge uniformly on \([0,1]\times[-1,1]\)
		to \(\Trans{\phi_{R,V}}(x,y)\).
	\end{enumerate}
\end{proposition}

The proof of \Cref{prop:wings} is deferred to \Cref{sec:wings}. One would want to construct the potential
\begin{equation}
	V_m^{\text{ideal}}(x,y) = 
	\begin{cases}
		V(x,y) \text{ if } (x,y)\in R\\
		m^2x+  g(m^{-1}(|x|-\pi/2),y) \text{ otherwise}.
	\end{cases}
\end{equation}
This potential, however, may not be convex, or continuous. The proof of \Cref{prop:wings} constructs the potentials \(\Wing{V}{g}\), which are convex approximations to the potentials \(V_m^{\text{ideal}}\).

The potentials \(\Wing{V}{g}\) go to \(+\infty\) in \(R_m \setminus R\). Therefore the measures \(\mu_{R_m,\Wing{V}{g}}\) converge strongly to the measure \(\mu_{V,\Omega}\). This implies that the first eigenfunction and eigenvalue of \((V_m, R_m)\) should converge as \(m\to \infty\) in \(L^2(R)\) to the first eigenfunction and eigenvalue of \((R,V)\). Completing this argument (in \Cref{sec:wings}) shows items \emph{1. -- 3}. of \Cref{prop:wings}.

The functions \(T_m:=\phi_{R_m, \Wing{V}{g}}(M x+\pi/2, \pi y)\) satisfy the PDE
\begin{equation}
	(1, \partial_y g(x,y)) \nabla T_m = \underbrace{m^{-1}(\partial_x g(x,y)\partial_x+\partial_{yy} + m^{-2}\partial_{yy}+\lambda_{R_m, V_m})}_{\text{Formally} \to 0}
\end{equation}
and therefore, once enough regularity is established (see \Cref{sec:wings}), one expects \(\lim_{m\to\infty}(1, \partial_y g(x,y)) \nabla T_m =0\). This means that any subsequential limit of \(T_m\) is constant on the flow lines of the vector field \((1, \partial_y g(x,y))\), which are curves of the form \((x,y(x))\), where \(y(x)\) satisfies the ODE \(y' = \partial_y g(x,y)\).  This property uniquely determines the limit to be \(\Trans{\phi_{R, V}}(x,y)\).

In our case we will consider \(g\) as follows:
\begin{lemma}\label{lem:flow_exists}
	One can choose \(g(x,y):[0,1]\times[-1,1]\to \R\) that is symmetric and convex in the second variable, and equal to zero in \([0,1]\times[-2/3,2,3]\) such that for all \(y\in [-1,1]\) the solution to the ODE
	\begin{equation}
	\begin{cases}
		{\gamma^{(y)}}'(x) = \partial_y g(x,\gamma^{(y)}(x))\\
		\gamma^{(y)}(1)=y
	\end{cases}
	\end{equation}
	satisfies \(\frac 23 \le |\gamma^{(y)}(0)|\le \frac 34 \).

	For \(q\) as in \Cref{def:hotspots_proposal}, and \(\beta\) as in \Cref{prop:perturbation} (so that \(\beta(\pi/2, y)-q(y)\) is equal to a constant \(c_0\)) the function \(\Trans q(x,y)\) is non-increasing in \(x\), and satisfies \(\Trans \beta(0,y) = q(y)+c_0\) and \(\Trans q(1,y) = q(\min(|y|, \frac 34))+c_0\). In particular \(\Trans \beta(0,y)-\Trans \beta(1,y)\) is a nonnegative continuous function that is not exactly zero.
\end{lemma}

\begin{proof}
	\textbf{Existence of \(g\):} Let \(g_0(x)\) be any symmetric convex smooth function that vanishes exactly on \([-2/3,2/3]\), and \(b(x)\) be a smooth function that vanishes on \([-1,1/2]\) and is positive on \((1/2,1]\). For \(M\) large enough, the function \(Mb(y)g_0(x)\) satisfies the hypothesis.

	\textbf{ \(\Trans \beta\) is non-increasing:}We want to show that \(\partial_x \Trans \beta \le 0\).  By symmetry, it suffices to consider \(x\ge 0\).  By construction of \(\Trans \beta\), 
	\begin{equation}
		\nabla \Trans  \beta \cdot (1, \partial_x g) =\partial_x \Trans  \beta + \partial_x g \partial_y \Trans \beta = 0,
	\end{equation}
	 and \( \partial_x g(x,y)\ge 0\) whenever \(y\ge 0\), and \(\partial_x g(x,y)= 0\) when \(|y|<\frac 2 3\). It suffices to show that \(g(x,y)\) is nondecreasing in \(y\) whenever \(y\ge \frac 23\). But \(\Trans \beta(x,y) = q(F(x,y))+c_0\), with \(F\) as in \Cref{def:transport}. If \(y>\frac 23\), then \(F(x,y)\ge 23\) and \(F\) is increasing in \(y\). Since \(q\) is nondecreasing in \(y\) for \(y\ge 23\), then \(g(x,y)\) is nondecreasing in \(y\).

	 The fact that \(\Trans \beta(0,y) = q(y)+c_0\) follows from the definition of \(\Trans \beta\). The fact that  \(\Trans q(1,y) = q(\min(|y|, \frac 34))+c_0\) follows from the fact that  \(\Trans q(1,y) = q(\gamma^{(y)}(0))+c_0\). The function \(\gamma^{(y)}(0)\) is equal to \(y\) if \(|y|<\frac 23\). Otherwise, \(|\gamma^{(y)}(0)|\in [2/3,3/4]\). Since \(q\) is constant in \([2/3,3/4],\) the result follows. 
\end{proof}

The flow associated to \(g\) arising from  \Cref{lem:flow_exists} and the function \(T_{R, \epsilon V_q, g}\), with \(q\) arising from \Cref{def:hotspots_proposal} look like
	\begin{center}
	\qquad
	\includegraphics[width=.15\textwidth]{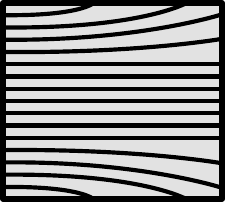}
	\qquad
	\includegraphics[width=.3\textwidth]{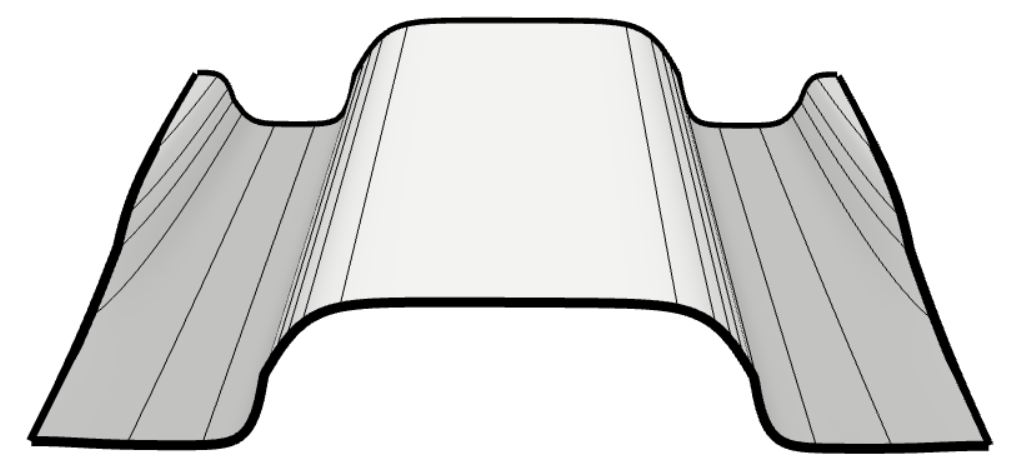}
	\end{center}
In the drawing of \(T_{R, \epsilon V, g}(x,y)\), the nearest side corresponds to \(x=1\), and the furthest to \(x=0\). Note how the level sets drawn on \(F\) correspond to the flow lines of the ODE on the left. The qualitative difference between \(T_{R, \epsilon V, g}(0,y)\) and \(T_{R, \epsilon V, g}(1,y)\) is that the later has \emph{lost} the high boundary spots.

\subsection{\notoc{Showing that \(h_{\epsilon,m}\) gives a counterexample.}}

\begin{definition}[Shortening notation]
	Fix any \(q\) satisfying \Cref{def:hotspots_proposal}, and \(V_q\), \(\beta\) arising from \Cref{prop:perturbation} from this \(q\). Fix any \(g\) satisfying \Cref{lem:flow_exists}. Let \(V_{\epsilon, m}:= \Wing{\epsilon V}{g}\), and by \(\phi_{\epsilon, m}:= \phi_{R_m, V_{\epsilon, M}}\) the ground eigenfunction of \((R_m, V_{\epsilon, m})\) with the sign chosen so that \(\phi_{\epsilon, m}>0\) for \(x>0\) (one can do that as long as \(\epsilon\) is small enough, and \(m\) is large enough, and the eigenfunctions converge to \(\sqrt{2/\pi}\sin(x)\)). We will denote by \(h_{\epsilon, m}\)  the heat extension \(h_{\epsilon, m}(x,y,t):=h_{R_m, V_{\epsilon, m}}(x,y,t)\). We will denote by \(\lambda_{\epsilon,m}:=\lambda_{R_m, V_{\epsilon, m}}\), and \(\lambda_\epsilon = \lambda_{R, \epsilon V_q} = \lim_{m\to \infty} \lambda_{\epsilon, m}\).
\end{definition}

\Cref{prop:heat_flow_spot} (and from it \Cref{thm:counterexample}) follows by showing that for \(\epsilon\) small enough, and \(V= V_{\epsilon, q}\) as in \Cref{prop:perturbation}, the inequality
\begin{equation}
	\lim_{m\to \infty} \max_{(x,y,t)\in R_m\times [0,1]}
	h_{\epsilon, m}>
	\lim_{m\to \infty} \max_{(x,y,t)\in \partial_h R_m}
	h_{\epsilon, m}
\end{equation}
holds. The definition below is key to estimating both limits:

\begin{definition}[Core limit, wings limit]
	 Let \(h_{\epsilon}^c(x,y,t):\R \times [-1,1]\times[0,1]\) be the unique continuous bounded solution to
	\begin{equation}
	\begin{cases}
		h^c_{\epsilon}(x,y,0) = \phi_{R,\epsilon V}\left(\frac{x}{\max(1,|x|)}, y, 0\right)\\
		\partial_t h^c_{\epsilon} = 
         			\frac{1}{8}(\Delta_{x,y}h^c_{V,g}+ \lambda_\epsilon h^c_{\epsilon})
         		& 
         		\text{ on } \R \times [-1,1]\times(0,1]
         	\\
         		\partial_{y} h^c_\epsilon = 0 
         		& 
         		\text{ on } \R\times \{-1,1\} \times (0,1].
	\end{cases}
	\end{equation}
	and 
	 \(h_\epsilon^w(x,y,t):[0,1]\times [-1,1]\times[0,1]\) be the unique continuous solution to 
	\begin{equation}
	\label{eq:heat_wings}
	\begin{cases}
		h^w_{\epsilon}(x,y,0) =\Trans{\phi_{R,\epsilon V}}\\
		\partial_t h^w_\epsilon = 
         			\frac{1}{8}(\partial_{yy} h^w_\epsilon+ \lambda_\epsilon h^w_\epsilon)
         		& 
         		\text{ on } [0,1] \times [-1,1]\times(0,1]
         	\\
         		\partial_{y} h^w_\epsilon = 0 
         		& 
         		\text{ on } [0,1]\times \{-1,1\} \times (0,1].
	\end{cases}
	\end{equation}

\end{definition}

\begin{proposition}
	\label{prop:wing_heat_extension}
	For each \(m\ge 0\), let \((x_m,y_m, t_m) \in R_m\times[0,1]\). Then
	\begin{enumerate}
	 	\item If the limit \(\lim_{m\to \infty} (x_m, y_m, t_m) =(x_\infty,y_\infty,t_\infty)\) exists and is finite, then
	\begin{equation}
	\label{eq:wing_limit_in}
	 	\lim_{m\to \infty} h_{\epsilon,m}(x_m, y_m, t_m) = h^c_{\epsilon}(x_\infty,y_\infty,t_\infty).
	 \end{equation} 
	 \item If \(\lim_{m\to \infty} x_m = \pm \infty\) and \(\lim_{m\to \infty} (x_m/m, y_m, t_m) =(x_\infty,y_\infty,t_\infty)\), then 
	\begin{equation}
		\label{eq:wing_limit_out}
	 	\lim_{m\to \infty} h_{\epsilon,m}(x_m, y_m, t_m) = \pm h^w_\epsilon(|x_\infty|,y_\infty,t_\infty).
	 \end{equation} 
	 \end{enumerate} 
\end{proposition}

\begin{proof}
	Consider the function \(\tilde h_{\epsilon,m}: \R \times [-1, 1] \times [0,1]\) to be the periodic extension in the \(x\) variable of \(h_{\epsilon,m}\). Since \(h_{\epsilon,m}\) has Neumann boundary conditions, the function \(\tilde h_{\epsilon,m}\) satisfies the PDE
	\begin{equation}
		\partial_t \tilde h_{\epsilon,m} = 
         			\frac{1}{8}(\Delta_{x,y}\tilde h_{\epsilon,m}+ \lambda_{\epsilon,m}\tilde h_{\epsilon,m})
	\end{equation}
	on \(R \times [-1-1]\times [0,1]\) with parabolic Neumann boundary conditions as well. The functions \(\tilde h_{\epsilon,m}(x,y,0)\) are uniformly bounded and converge in \(L^\infty_{loc}(\R \times [-1, 1])\) to \(\phi_{V,g}\left(\frac{x}{\max(1,|x|)}, y, 0\right)\) by \Cref{prop:wings}, from which~\eqref{eq:wing_limit_in} follows.

	The case~\eqref{eq:wing_limit_out} follows from the fact that \((x,y) \mapsto \tilde h_{\epsilon,m}(x-x_m,y,0)\) converges locally uniformly in \(\mathbb R \times [-1,1]\) to \(\pm h^w(|x_\infty|, y,0)\) (with positive sign if \(x_m \to +\infty\), and negative sign otherwise). Since the limit is constant in the \(x\) variable, the heat flow becomes a single-variable heat flow in the \(y\) variable, as in~\eqref{eq:heat_wings}.
\end{proof}

From the sign conventions that \(\phi_{R} = +\sqrt{2/\pi}\sin(x)\), that \(\phi_{R,\epsilon V} \to \phi_R\) as \(\epsilon \to 0\) (\Cref{prop:perturbation}) implies that for \(\epsilon\) small enough, \(\phi_{R,\epsilon V}(\pi/2,y)>0\), and therefore \(\Trans{\phi_{R,\epsilon V}(\pi/2,y)}>0\). In particular, in our sign convention, \(h^w_{\epsilon}>0\) for all \(\epsilon\) small enough.

\begin{corollary}
	Assume \(\epsilon\) is small enough (depending on \(g,q,V\)). Let \(R^w := [0,1]\times [-1,1]\times[0,1])\) (the domain of \(h^w_{V,g}\)), and \(R^c := \R \times [-1,1]\times[0,1])\) (the domain of \(h^c_{\epsilon}\)). 

	Let \(D^w\) be the subset of \(R^w\) where \(x=1\), or \(y\pm1\), or \(t=0\), and \(D^c\) the subset of \(R^c\) where either \(t=0\) or \(y\pm1\).

	Then
	\begin{equation}
	\label{eq:limit_interior}
		\lim_{m\to \infty}\max_{((x,y),t)\in R_{m}\times [0,1]} h_{\epsilon,m}(x,y,t) = \max\left(\max_{((x,y),t)\in R^w} h_\epsilon^w(x,y,t), \max_{((x,y),t)\in R^c} h_\epsilon^c(x,y,t)\right)
	\end{equation}
	and
	\begin{equation}
	\label{eq:limit_boundary}
		\lim_{m\to \infty} \max_{((x,y),t)\in \partial_h R_m} h_{\epsilon,m}(x,y,t) = \max\left(\max_{
		(x,y,t)\in D^w} h_\epsilon^w(x,y,t), \max_{(x,y,t)\in D^c} h_\epsilon^c(x,y,t)\right).
	\end{equation}
\end{corollary}
\begin{proof}
	Every sequence of points \((x_m, y_m, t_m)\) \(R_m\times[0,1]\) has a subsequence satisfying either case 1. or 2. in \Cref{prop:wing_heat_extension}. This allows one to characterize the limit of the maxima of \(h_{\epsilon,m}\) as \(m\to \infty\). Taking \((x_m, y_m, t_m)\) to be the extremizers of \(h_{\epsilon, m}\) in either \(\partial_h R_m\) or \(R_m \times[0,1]\) gives the result.
\end{proof}

To combine \Cref{prop:wing_heat_extension} with \Cref{prop:perturbation} we show that the maximum value of \(h_{\epsilon,m}\) in the boundary must be (up to an \(o(\epsilon)\)) attained in the wings.

\begin{proposition}[Maximum is in the wings]
\label{prop:wing_bdy_larger}
	Let \(\beta(x,y)\) be a smooth perturbation, and \(V_{\epsilon}\) as in \Cref{prop:perturbation}. Fix a smooth convex function \(g\) as in \Cref{prop:wings}. Then
	\begin{equation}
		\max_{
		(x,y,t)\in D^c} h_{\epsilon}^c(x,y,t)\le  \max_{(x,y,t)\in D^w} h_{\epsilon}^w(x,y,t)+o(\epsilon).
	\end{equation}
\end{proposition}

\begin{proof}
	By the maximum principle applied to the heat flow, it suffices to show that
	\begin{equation}
		\max_{(x,y)\in R} h^c_\epsilon (x,y,0) - h_{\epsilon}^w(0,y,0) = o(\epsilon).
	\end{equation}
	By expanding definitions of \(h_{\epsilon}^c, h_\epsilon^w\), it suffices to show
	\begin{equation}
		\max_{(x,y)\in R} \underbrace{\sqrt{2/\pi}\sin(x) +\beta(x,y)}_{h_{\epsilon}^c(x,y,0)+o(\epsilon)} - \underbrace{\sqrt{2/\pi}  -\beta(\pi/2, y)}_{h_{\epsilon}^w(0,y,0)+o(\epsilon)} = o(\epsilon).
	\end{equation}
	this in turn would follow from
	\begin{equation}
		\max_{(x,y)\in R} \beta(x,y)-\beta(\pi/2,y) - \epsilon^{-1} \sqrt{2/\pi}(1-\sin(x)) =o(1).
	\end{equation}
	Since \(\beta\) is Lipschitz, there a constant \(C\) such that \(\beta(x,y)- \beta(\pi/2,y) +o(1) \le C (\pi/2-x)\), and therefore 
	\begin{equation}
		\max_{(x,y)\in R} \beta(x,y)-\beta(\pi/2,y) - \epsilon^{-1} 2/\pi(1-\sin(x)) =o(1).
	\end{equation}
\end{proof}

This allows one to drop \(h_{\epsilon}^c(x,y,t)\), and study only \(h_{\epsilon}^w(x,y,t)\) in order to show \Cref{thm:counterexample}.
\begin{proposition}
\label{prop:final_approximation}
	The functions
	\(p_{\epsilon}(x,y,t):=\epsilon^{-1}(h^w_{\epsilon}(x,y,t)\exp(-\lambda_\epsilon/8 \cdot t)-\sqrt{2/\pi})\) converge uniformly as \(\epsilon\to 0\) in \([0,1]\times [-1,1]\times[0,1]\) to a function \(p_0\) satisfying
	\begin{equation}
		\begin{cases}
			p_0(x,y,0)= \Trans{\beta}\\
			\partial_t p_0 = \frac 1 8 \partial_{yy} p_0 \text{ in } [0,1]\times [-1,1]\times[0,1]\\
			\partial_y p_0 = 0 \text{ in } [0,1]\times \{-1,1\}\times[0,1]
		\end{cases}
	\end{equation}
	Where \(\beta\) is as in \Cref{prop:perturbation}, and satisfies \(\beta(\pi/2,y) = q(y)+c\) for some fixed constant \(c\).
\end{proposition}

\begin{proof}
	The functions \(h^w_{\epsilon}(x,y,t)\exp(-\lambda_\epsilon t)\) satisfy the heat equation in the \(y\) variable. By the maximum principle, it suffices to show uniform convergence of the initial data (\(t=0\)).

	The map \(f\mapsto \Trans f\) is linear bounded (with norm 1) from \(L^\infty(\{\pi/2\}\times[-1,1])\) to \(L^{\infty}([0,1]\times[-1,1])\). In particular, it suffices to show that \(\epsilon^{-1}(\phi_{R, \epsilon}(\pi/2,y)-\sqrt{2})\) converges uniformly in \([-1,1]\) to \(q(y)\). This is a the content of~\eqref{eq:approx_limit}.

\end{proof}

We now proceed to the proof of \Cref{thm:counterexample}.

\begin{proof}[Proof of \Cref{thm:counterexample}.]

	By \Cref{prop:wing_bdy_larger} it is enough to show that there exists a small constant \(c\) such that for all \(\epsilon\) small enough, 

	\begin{equation}
	\label{eq:final_goal}
		\max_{(x,y,t)\in R^w} h_{\epsilon}^w(x,y,t)-
		\max_{(x,y,t)\in D^w} h_{\epsilon}^w(x,y,t) > c\epsilon.
	\end{equation}

	For \(0\le t\le 1/2\), one has the estimate  \(h_\epsilon^w(x,y,t) \le \exp(\lambda_0/16) + O(\epsilon)\) uniformly in \(x,y\). On the other hand, \(h_\epsilon^w(x,y,1) = \exp(\lambda_0/8) + O(\epsilon)\) (again uniformly). In particular, both maxima in the expression above are taken for \(t>1/2\). This rules out that the maximum on \(D^w\) is taken in the part where \(t=0\). Thus we have to show two estimates, one for the case \(\max_{(x,y,t)\in D^w} h_{\epsilon}^w(x,y,t)\) is attained at \(x=1\), and another one if it is attained at \(y = \pm 1\). These estimates are
	\begin{equation}
	\begin{cases}
		\max_{\substack{y\in [-1,1]\\t\in [1/2,1]}} \epsilon^{-1}(h_{\epsilon}^w(0,y,t)-
		 h_{\epsilon}^w(1,y,t))  \ge c \\
		\max_{\substack{x\in [0,1]\\t\in [1/2,1]}} \epsilon^{-1}(h_{\epsilon}^w(0,0,t)-
		 h_{\epsilon}^w(x,\pm 1,t))  \ge c 
	\end{cases}
	\end{equation}
	and together imply~\eqref{eq:final_goal}. By \Cref{prop:final_approximation}, one can send \(\epsilon\to 0\), where is enough to show

	\begin{equation}
	\begin{cases}
		\max_{\substack{y\in [-1,1]\\t\in [1/2,1]}} 
		p_0(0,y,t)-
		p_0(1,y,t))  \ge c \\
		\max_{\substack{x\in [0,1]\\t\in [1/2,1]}} p_0(0,0,t)-
		 p_0(x,\pm 1,t))  \ge c 
	\end{cases}
	\end{equation}

	\textbf{First estimate:} The function \(p(0,y,t)-p(1,y,t)\) solves the Neumann heat equation with initial data \(\Trans \beta (0,y)- \Trans \beta(1,y)\). By \Cref{lem:flow_exists}, this initial data is continuous, nonzero, and nonnegative, and therefore by the strong maximum principle \(p(0,y,t)-p(1,y,t)\) is bounded from below by a constant in \([-1, 1]\times[1/2, 1]\).

	\textbf{Second estimate:} By the weak maximum principle, and the fact that \(\Trans \beta\) is decreasing in the \(x\) variable, for any \(x\in [0,1]\) one has \(p(x,y,t)\le p(x,0,t)\), and thus the second inequality follows by showing that for any \(t\in [1/2,1]\)
	\begin{equation}
	\label{eq:2nd_bdy_bound}
		p(0,0,t) - p(0,\pm1,t) \ge c,
	\end{equation}
	but 
	\(p_t(0,y,t)\) is (up to an additive constant) precisely the heat flow of \(q(y)\). The construction of \(q(y)\) (the condition \(H_q(0,t)-H_q(1,t)>0\) in \Cref{def:hotspots_proposal}) is equivalent to 
	\(p(0,0,t) - p(0,\pm 1,t)>0\) in \([0,1]\). By continuity, there is a constant such that~\eqref{eq:2nd_bdy_bound} holds, finishing the proof.
\end{proof}

%% file: Parts/basic_cylinder_properties.tex
The goal of this section is to show the following parts of \Cref{thm:cylindrical_convergence}:

\begin{proposition*}[Copy of the first parts of \Cref{thm:cylindrical_convergence}]
    Let \((\Omega, V)\) be a smooth convex pair with \(V\ge 0\) a spectral gap, and \(\lambda_{\Omega,V}\le 2\). Let \(\Omega_d:=\Barrel{\Omega}{V}\), then:
    
 \begin{enumerate}[label ={(C.\arabic*)}]
     \item The ground eigenvalues \(\lambda_{\Omega_d}\) converge to \(\lambda_{\Omega,V}\) as \(d\to \infty\). For \(d\) large enough, the sets \(\Omega_d\) have a spectral gap. \
     \item For \(d\ge 4\|V\|_\infty\), the ground eigenfunction  \(\phi_{\Omega_d}(x,w)\) of \((\Omega_d,0)\) is radial in the second coordinate, and thus can be written in the form  \(\phi_{\Omega_d}(x,w)=: \psi_d(x,1-(|w|/R(d))^2)\).
\end{enumerate} 
\end{proposition*}

The key geometric aspects of \(\Omega_d\) (which are the main design considerations for \(\Omega_d\)) are the following two elementary Lemmas:

\begin{lemma}[The measures \(\1_{\Omega_d}\) converge to \(\exp(-V) \1_\Omega\)]\label{lem:measures_converge}
    The following asymptotic equality holds
    \begin{equation}
        \label{eq:measure_approx}
        \frac{|\{(x,w) \in \Omega_d, x= x_0\}|}{|B_{\sqrt{d}/2}(\R^{d+1}))|} = e^{-V(x_0)}(1+O_{\|V\|_\infty}(d^{-1})).
    \end{equation}

    In particular, if, as in \Cref{def:barrel}, \(\mu_{\Omega_d}:= \frac 1 {|\Omega_d|} \1_{\Omega_d}\) is the normalized uniform measure in \(\Omega_d\) (and \(\mu_{\Omega,V}\) is the weighted probability measure with weight \(\frac{1}{Z_{\Omega,V}}\1_{\P,ega}\exp(-V)\)), setting \(\pi^1(x,w) = x\), then
    \begin{equation}\pi^1_*\mu_{\Omega_d}\to\mu_{\Omega,V}\end{equation}strongly.
\end{lemma}

\begin{proof} 
    \begin{align}
        \frac{H^{d+1}(\{(x,w) \in \Omega_d, x= x_0)}{H^{d+1}(B_{\sqrt{d}}(\R^{d+1}))}
        = \left(1+\frac{V(x_0)}{d}\right)^{d+1} =  e^{-V(x_0)}(1+O_{\|V\|_\infty}(d^{-1})).
    \end{align}
\end{proof}

\begin{lemma}[The balls \(B_{R_d(x)}\) are spectrally small]\label{lem:balls_spectrally_small}  For \(d\ge 4 \| V(x)\|\) and all \(x\in \Omega\), and any \(r\le \frac 1 2 \left(\sqrt{d}-\frac{V(x)}{\sqrt{d}}\right)\) the the first nonzero Neumann eigenvalue of the ball \(B_r(\R^{d+1})\) is at least \(2.5\).

\end{lemma}

\begin{proof}

    If \(d\ge 4 \|v\|\), then \(r\le \frac 1 2 \sqrt{d}\left(1+\frac{1}4\right) = \frac{5}{8}\sqrt{d}\). Since \(\lambda_{c\Omega} = c^{-2}\lambda_{\Omega}\), it suffices to show that \(\lambda_{B_1(\R^{d+1})}\ge d\) (See \Cref{lem:poincare_ball}). Then for \(r\le \sqrt{d}\frac 1 2 \left(1+\frac{1}4\right)\),
    \begin{equation}
        \lambda_{B_r(\R^d)} \ge \left(\frac{5\sqrt{d}}{8}\right)^{-2}\! \!\!\cdot d = 2.56.
    \end{equation} 
\end{proof}

\begin{lemma}
    Let \(f\) be \(L^2\) function \(\Omega_d\) with mean zero. Let \(f(x,w)= f_0 + \bar f(x) +\tilde f(x,w)\)  where 
    \begin{enumerate}
        \item The function \(\tilde f\) has mean zero on the balls given by \(\{(x,w) \in \Omega_d, x=x_0\}\).
        \item The function \(\bar f(x)\) is constant in \(w\) and has mean zero in \(\Omega\) with respect to \(\mu_{\Omega,V} = \frac{1}{Z_{\Omega,V}}\exp(-V) dx\).
        \item The remainder \(f_0\) is a constant.
    \end{enumerate}
    Then
    \begin{align}
        f_0 = O_{\Omega,V}(d^{-1}) \|\bar f\|_{L^2(\mu_{\Omega_d})}.
    \end{align}
    \begin{equation}
    \label{eq:almost_orth_l2}
        \|f\|_{L^2(\mu_{\Omega_d})}^2 = (1+O_{\Omega,V}(d^{-1}))(\|\tilde f\|_{L^2(\mu_{\Omega_d})}^2  + \|\bar f\|_{L^2( \mu_{\Omega,V})}^2 )
    \end{equation}
    \begin{equation}
    \label{eq:almost_orth_h1}
        \|\nabla f\|_{L^2(\mu_{\Omega_d})}^2 = (1+O_{\Omega,V}(d^{-1/2}))(\|\nabla_w \tilde f\|_{L^2(\mu_{\Omega_d})}^2  + \|\nabla_x\bar f\|_{L^2(\Omega, \mu_{\Omega,V})}^2 )
    \end{equation}

\end{lemma}

\begin{proof}

    By \Cref{lem:measures_converge}, for any function \(f:\Omega\to \mathbb R\), the quantities
    \begin{equation}
        \|f\|_{L^2(\mu_{\Omega_d})}, \quad
        \|f\|_{L^2(\mu_{\Omega,V})},
    \end{equation}
    where in the first one \(f\) is understood to be constant in the \(w\) variable, differ by a factor of at most \((1+O_{\Omega,V}(d^{-1}))\). Thus, these norms are interchangeable.

    The function \(\bar f(x)\) must also have mean zero on \(\Omega_d\). Integrating the \(w\) variable, \(\bar f(x)\) has mean zero in \(\Omega\) with respect to the weight \( \left(1+\frac{V(x)}{d}\right)^{d+1}\). On the other hand
    \begin{align}
        f_0 =& \frac 1 {Z_{\Omega,V}} \int_\Omega \bar f(x) \exp(-V(x)) dx
        \\= &
        \frac 1 {Z_{\Omega,V}} \int_\Omega \bar f(x)\left( \exp(-V(x)) - \left(1+\frac{V(x)}{d}\right)^{d+1} \right) dx
        \\=& O_{\Omega,V}(d^{-1}) \|\bar f\|_{L^2(\mu_\Omega)}
        \\\le & e^{\|V\|_{\infty}}O_{\Omega,V}(d^{-1}) \|\bar f\|_{L^2(\mu_{\Omega,V})}.
    \end{align}

    The functions \(\bar f, \tilde f\) are orthogonal in \(L^2(\mu_{\Omega_d})\), and therefore 
    \begin{equation}
        \|f\|_{L^2(\mu_{\Omega_d})}^2 = (1+O_{\Omega,V}(d^{-1}))(\|\tilde f\|_{L^2(\mu_{\Omega_d})}^2  + \|\bar f\|_{L^2(\mu_{\Omega_d})}^2 )
    \end{equation}

    An explicit computation shows that
    \begin{align}
        \grad_x \bar f = 
        \fint_{B_{R_d(x)}} 
        \grad_x f(x,w) +\frac 1d (1+d^{-1} V(x))^{-1}\grad_x V (w\cdot \grad_w f) dw.
    \end{align}
    Using the elementary inequality \((a+\delta b)^2 \le a^2+ 2\delta (a^2+b^2)\) this gives, for \(d\) large enough,
    \begin{equation}
        \int_{\mu_{\Omega_d}} |\grad_x \bar f|^2 dx dw \le \int_{\mu_{\Omega_d}}  |\grad_x f|^2 + \frac{8\|\grad V\|_{L^{\infty}}}{\sqrt{d}} |\grad f|^2 dx dw.
    \end{equation}
    Therefore, using that \(\nabla_w f = \nabla_w \tilde f\),
    \begin{align}
        \|\nabla f\|_{L^2(\mu_{\Omega_d})}^2 =&
        \|\nabla_w f\|_{L^2(\mu_{\Omega_d})}^2 +
        \|\nabla_x f\|_{L^2(\mu_{\Omega_d})}^2
        \\=&
         (1+O_{\Omega,V}(d^{-1/2}))(\|\nabla_w \tilde f\|_{L^2(\mu_{\Omega_d})}^2  + \|\nabla_x\bar f\|_{L^2(\mu_{\Omega_d})}^2 )
    \end{align}
\end{proof}

\begin{proof}[Proof of \Cref{prop:cylindrical_convergence}~\ref{item:eigenvalues_converge}]

    A direct test on the Rayleigh quotient using \Cref{lem:measures_converge} shows that \(\lambda_{\Omega_d} \le \lambda_{\Omega,V}(1+O(d^{-1}))\):
    \begin{align}
    \lambda_{\Omega_d} = 
    &
    \min_{\substack{
        f\in H^{1}(\Omega_d)\\ \int_{\Omega_d} f(x,w) \dd x \dd w = 0
        }} 
    \frac
        {\int_{\Omega_d} |\grad f(x,w)|^2 \dd x \dd w}
        {\int_{\Omega_d} | f(x,w)|^2 \dd x \dd w}
    \\\le&
    \min_{\substack{
        f\in H^{1}(\Omega)\\ \int_{\Omega_d} f(x) \dd x \dd w = 0
        }} 
    \frac
        {\int_{\Omega_d} |\grad f(x)|^2 \dd x \dd w}
        {\int_{\Omega_d} | f(x)|^2 \dd x \dd w}
    \\=&
    (1+O(d^{-1}))
    \min_{\substack{
        f\in H^{1}(\Omega)\\ \int_{\Omega} f(x) \exp^{-V(x)}\dd x = 0
        }} 
    \frac
        {\int_{\Omega_d} |\grad f(x)|^2 e^{-V(x)} \dd x}
        {\int_{\Omega_d} | f(x)|^2 e^{-V(x)} \dd x}.
    \end{align}

    Let \(\phi_{\Omega_d}\) be the first eigenfunction of \(\Omega_d\), of eigenvalue \(\lambda_{\Omega_d} \le \lambda_{\Omega,V}+O(d^{-1})\), normalized so that \(\|\phi_{\Omega_d}\|_{L^2(\Omega_d)}=1\). Then, combining~\eqref{eq:almost_orth_l2} and~\eqref{eq:almost_orth_h1}:
    \begin{align}
    \lambda_{\Omega_d} = &
    \|\nabla\phi_{\Omega_d}\|_{L^2(\Omega_d)}
    \\=& 
    (1+O_{\Omega,V}(d^{-1/2}))
    (\|\nabla_w \tilde \phi_{\Omega_d}\|_{L^2(\Omega_d)}
    +\|\nabla_x \overline{\phi_{\Omega_d}}\|_{L^2(\Omega, \exp(-V))}
    )
    \\\ge&
        (1+O_{\Omega,V}(d^{-1/2}))
    (2.5\|\tilde \phi_{\Omega_d}\|_{L^2(\Omega_d)}
    +\lambda_{\Omega,V}\|\overline{\phi_{\Omega_d}}\|_{L^2(\Omega, \exp(-V))}
    )
    \\\ge&
    (1+O_{\Omega,V}(d^{-1/2}))
    \lambda_{\Omega,V}\|\phi_{\Omega_d}\|_{L^2(\Omega_d)} + (2.5-\lambda_{\Omega_d})\|\tilde \phi_{\Omega_d}\|_{L^2(\Omega_d)}
    \end{align}

    This shows that \(\lim_{d\to \infty} \lambda_{\Omega_d} = \lambda_{\Omega, V}\).  Let \(\lambda^{(2)}_{\Omega,V}\) be the second eigenvalue of \((\Omega,V)\), counted with multiplicity. 
    Let \(f\) be any function in \(\Omega_d\) satisfying both \(\int_{\Omega_d} f(x,w) dx dw = 0\), and \(\int_{\Omega} \bar f \phi_{\Omega,V}\exp(-V) dx = 0\). The same argument shows
    \begin{align}
    \|\nabla f\|_{L^2(\Omega_d)}
    & 
    (1+O_{\Omega,V}(d^{-1/2}))
    (\|\nabla_w \tilde f\|_{L^2(\Omega_d)}
    +\|\nabla_x \overline{f}\|_{L^2(\Omega, \exp(-V))}
    )
    \\\ge&
        (1+O_{\Omega,V}(d^{-1/2}))
    (2.5\|\tilde f\|_{L^2(\Omega_d)}
    +\lambda^{(2)}_{\Omega,V}\|\overline{f}\|_{L^2(\Omega, \exp(-V))}
    )
    \end{align}
    Since \((\Omega,V)\) has a spectral gap, \(\lambda^{(2)}_{\Omega,V}>\lambda_{\Omega,V}\). By the variational characterization of eigenvalues \(\Omega_d\) must have a spectral gap for large \(d\).
\end{proof}

\begin{lemma}%
\label{lem:eigenfunctions_radial}
    Let  \(d\ge 4\|V(x)\|\). Let \(f(x,w)\) be a Neumann eigenfunction of \(\Omega_d\) of eigenvalue \(\lambda<2.5\). Then \(f\) is radial in the second coordinate.
\end{lemma}

\begin{proof}
    Assume \(f\) was not radial in the second coordinate. Let \(\bar f(x,r) := \fint_{\mathbb S^{d}} f(x, r\vec e) d \vec e\). Then \(\tilde f(x,w):= f(x,w)- \bar f(x,w)\) is a nonzero Neumann eigenfunction of the Laplacian with the same eigenvalue. For each fixed \(x_0\), the function \(\tilde f(x_0,w)\) has mean zero on the ball \(\{(x_0,w), |w|\le  \rho_V(x))\}\). By \Cref{lem:balls_spectrally_small}, 
   \begin{equation}
        \int_{|w|\le  \rho_V(x)} |\nabla_w f(x_0,w)|^2 dw \ge 2.5
        \int_{|w|\le  \rho_V(x)} | f(x_0,w)|^2 dw 
   \end{equation}
   Therefore
    \begin{equation}
        \lambda = \frac{\int_{\Omega} \int_{|w|\le  \rho_V(x)} |\nabla f|^2 dw dx }{\int_{\Omega} \int_{|w|\le  \rho_V(x)} |f|^2 dw dx}
        \ge  
        \frac{\int_{\Omega}2.5  \int_{|w|\le  \rho_V(x)} | f|^2 dw dx }{\int_{\Omega} \int_{|w|\le  \rho_V(x)} |f|^2 dw dx} = 2.5
    \end{equation}
\end{proof}

\begin{proof}[Proof of \Cref{prop:cylindrical_convergence}~\ref{item:eigenfunctions_radial}]
    If \(\lambda_{\Omega_d}<2.5\) (which happens for all \(d\) large enough), \(\phi_d\) has to be radial. Recalling the definition of \(S_d(\Omega,V)\) in \cref{def:barrel}, there is a function \(\phi^{S}_d:S_d(\Omega, V)\to \R\) such that
\begin{equation}
\label{eq:phi_s_def}
\phi_d(x,w)= \phi^S_d(x,\sqrt{d}/2-|w|),
\end{equation}
 where the function \(\phi^S_d\) is the lowest (non-constant, Neumann) eigenfunction of \(L^s_{\Omega,V}:=-\Delta+\grad V_s \cdot \grad \) in \(S_d(\Omega, V)\) with eigenvalue precisely \(\lambda_{\Omega_d}\).  The function \(\psi_d\) in~\ref{item:eigenfunctions_radial} is just a change of variables on \(\phi^r_d\).
\end{proof}

%% file: Parts/boundary_l2.tex
The goal of this section is to show the following proposition:

\begin{proposition*}[Part of \Cref{thm:cylindrical_convergence}, with changes underlined.]
    Let \((\Omega, V)\) be a convex pair with a spectral gap, with \(\lambda_{\Omega, V}\le 2\), and let \(\Omega_d:=B_d(\Omega, V)\), then:
     \begin{enumerate}
     \item[{(C.3)}] {If  \underline{a subsequence of} \(\psi_{d}(x,t)\) converges  weakly in \(L^2(\Omega \times[0,1])\) to a function \(\psi_\infty\), then \(\psi_\infty\) satisfies:}
     \begin{enumerate}
        \item[{(a)}] {\(\psi_\infty(x,0) = \phi_{\Omega, V}(x)\)} 
    \end{enumerate}
    \end{enumerate}
\end{proposition*}

In the process, we will show that \(\phi^S_d\) are uniformly bounded and equicontinuous, which will show that

\begin{proposition}%
\label{prop:equicontinuity}
    The functions \(\phi^{S}_d(x,s)\) (as defined in~\eqref{eq:phi_s_def}) are equicontinuous, and therefore 
    \begin{enumerate}
        \item For any \(C>0\), as \(d\to \infty\) the functions \(\phi^{S}_d(x,s)\) converge uniformly in the set \(\{(x,s)\in S_d(\Omega,V), s\le C\}\) to the function \(\phi_{\Omega,V}(x)\).
        \item The functions \(\psi_d(x,t)\) are equicontinuous in the \(x\) variable and converge  uniformly to \(\phi_{\Omega,V}(x)\) in the sets \(\Omega \times [0,d^{-1/2}]\).
    \end{enumerate}
\end{proposition}

By uniqueness of solutions to the heat equation, once (C.3) is proven for all subsequences, the sub-sequential limit becomes a limit. The proof uses the variational characterization of eigenfunctions, with respect to the measures \(\sqrt{d}\exp(-\sqrt{d}W^s_d(x,r))\). These measures become more and more singular. To extract a limit one needs to first show uniform boundedness of the minimizers. To do so, we need the following estimates on the slice measures \(\mu_{\Slice \Omega V,\sqrt{d}W_{d}}\).

\begin{lemma}[Estimates on the slice measure]%
\label{lem:slice_estimates}
The normalization constants 
\begin{equation}
    \tilde Z_d:= \sqrt{d} \int_{S_d(\Omega,V)} \exp(-\sqrt{d}W^s_d(x,r)) dx dr
\end{equation} converge to \(Z_{\Omega,V} :=\int_\Omega\exp(-V) dx\). Moreover the bounds
\begin{equation}
\label{eq:radon-nykodim-upper}  
    \max_{(x,r)\in \Slice \Omega V} \frac{\exp(-\sqrt{d}W^s_d(x,r))}{\exp(- \sqrt{d} \cdot 2r) } \le 1+o_{\Omega,V}(1)
\end{equation}
and
\begin{equation}
\label{eq:radon-nykodim-lower}  
    \min_{\substack{(x,r)\in \Slice \Omega V\\|r|<d^{-1/4}}} \frac{\exp(-\sqrt{d}W^s_d(x,r))}{\exp(- \sqrt{d} \cdot 2r) } \ge 1+o_{\Omega,V}(1)
\end{equation}
hold as \(d\to \infty\) and for all \((x,r)\in S_d(V,\Omega)\).
\end{lemma}

\begin{proof}
    The potentials \(W_d^s(x,r) = -\sqrt{d} \log(1-2/\sqrt{d}x)\) are decreasing in \(d\) and converge to the function \(f(x) = 2 x\). Therefore
    \begin{align}
        \lim_{d\to \infty} \tilde Z_d 
        =& 
        \lim_{d\to \infty} 
        \int_{\Omega} \int_{V(x) /2\sqrt{d}}^{\infty} \sqrt{d}e^{-\sqrt{d}W^s_d} \1_{\Slice{\Omega}{V}} dx dr
        \\=&
        \int_{\Omega}  \int_{V(x)/2}^{\infty} e^{-2x} dx dr  = Z_{\Omega, V}.
    \end{align}

    The inequality~\eqref{eq:radon-nykodim-upper} follows from the fact that \(W^s_d(r)\ge 2x\) for \(r>0\).
    If \((x,r)\in S_d\), then \(r>-O(d^{-1.2})\). Then for \(r<0\) the result follows from a Taylor expansion of \(W_d^s\):
    \begin{equation}
        -\sqrt{d}(W^s_d(r)-2r) = \sqrt{d}\left(\sqrt{d}\left(\frac{2r}{\sqrt{d}}+O\left(\frac{2r}{\sqrt{d}}\right)^2\right) - 2r\right) = O(r^2)
    \end{equation}

    The inequality~\eqref{eq:radon-nykodim-upper} follows from the same Taylor expansion.
\end{proof}

\subsubsection{Uniform Harnack estimates}

Let \(\Omega\subseteq \R^d\) be a convex, bounded set with nonempty interior, \(V:\Omega \to \R\) be a Lipschitz convex function, and \(P_t\) be the semigroup generated by \(\Delta -\nabla V \cdot \nabla\) with Neumann boundary conditions in \(\partial \Omega\). That is to say, \(P_t f\) solves the equation
\begin{equation}
    \partial_t (P_t f)(x) = (\Delta -\nabla V \cdot \nabla)(P_t f)
\end{equation}
with Neumann boundary data and \(P_0 f = f\). The Li-Yau inequality states that \(P_t(f)\), which is defined for any \(f\in L^2(\Omega, \exp(-V))\),  has the following property:
\begin{theorem}
    [Dimension Free Harnack inequality\footnote{\Cref{thm:li-yau_harnack} and \Cref{thm:non-explosiveness} are often presented for \(C^2\) convex sets and \(C^2\) potentials. The constants do not depend on the regularity of \(\Omega, V\) and they can be extended to much lower regularity settings by standard approximation procedures. See \Cref{lem:fokker_plank_stable} for a simple self-contained version of this fact in the generality needed for the current proof.}, \(\kappa = 0\). Originally by~\cite{li1986parabolic}, see the survey {\cite[Theorem 2.1]{wang2006dimension}}]%
\label{thm:li-yau_harnack}
     Let \(\Omega\) be a convex domain, and \(V\) a Lipschitz convex potential in \(\Omega\).
    For any \(\alpha > 1, t>0\) and \(f\in L^2(\Omega, \exp(-V))\) one has
    \begin{equation}
        (P_t^{(\Omega,V)}|f|)^\alpha(x) \le(P_t^{(\Omega,V)} |f|^\alpha)(y)\exp\left(\frac 1 {4t} \frac{\alpha}{\alpha-1}|x-y|^2 \right).
    \end{equation}
\end{theorem}
This inequality, together with the non-explosiveness of the flow of \(P_t^{(\Omega, V)}\), will guarantee the convergence of \(\psi_n\) to \(\psi\).
\begin{proposition}[Non-explosiveness]%
\label{thm:non-explosiveness}
    Let \(\eta_\phi(\delta):=\sup_{\substack{x,y}\in \Omega, \  |x-y|<\delta} |f(x)-f(y)|\), the modulus of continuity of \(f\). Then for any \(t>0\), \(\delta>0\), one has
    \(\eta_{P_t f}(\delta) \le \eta_{f}(\delta) \).
\end{proposition}
\Cref{thm:non-explosiveness} is a stronger version of the usual non-explosiveness, and follows from the same Brownian coupling proof (see the proof in \Cref{sec:contractive_flows}) as the more standard
\begin{corollary}
    For any \(t>0\) one has
    \(\|\nabla P_t f\|_\infty \le \|\nabla f \|_\infty\).
\end{corollary}

As a direct consequence of the Li-Yau inequality, following ~\cite{rockner2003supercontractivity}, we obtain the well known result:

\begin{lemma}[Ultracontractivity]%
\label{lem:ultracontractivity}
    Let \((\Omega,V)\) be a convex pair with measure \(\mu_{\Omega,V}:=Z_{\Omega,V} \exp(-V) \1_\Omega\) and heat semigroup \((P_t)_{t\ge 0}\), for which the Li-Yau inequality holds. Assume that there exists a \(t\ge 0\) and \(C>0\) such that
    \begin{equation}
        \|P_{t/2}(\exp(2|x|^2/t))\|_{\infty} \le M \text{, and }
        \int_{\Omega} \exp(-2|x|^2/t) d\mu \ge M^{-1}
    \end{equation}
    then \(\|P_t f\|_{L^\infty(\Omega)} \le M^2 \|f\|_{L^2(\mu))}\).
    In particular, for any \(t>0\),
    \(\|P_t f\|_\infty \le C_{t, \operatorname{diam}(\Omega)} \|f\|_{L^2(\Omega, \mu_{\Omega,V})}\).
\end{lemma}
\begin{corollary}
    If, under the hypothesis of \Cref{lem:ultracontractivity}, \(f\) is a Neumann eigenfunction of \(-\Delta + \nabla V\cdot\nabla\) of eigenvalue \(\lambda\) then \(\|f\|_{L^\infty}\le M^2e^{\lambda t}\|f\|_{L^2(\Omega,\mu_V)}\).
\end{corollary}
\Cref{lem:ultracontractivity} is proven as \Cref{lem:ultracontractivity_proof} in the appendix, following the usual proof, explicitly tracking the constants. We will apply this result to \(P^{(d)}_t:=P^{(S_d,(\Omega, V),W^s_d(\Omega,V))}_t\), the (Neumann) heat semigroup in \((S_d(\Omega,V), \) by \(-\Delta + \sqrt d W^s_d(r)\).

\begin{lemma}%
\label{lem:li_yau_holds}
    The convex pair \((S_d, \sqrt d W_d)\) satisfies the Li-Yau Harnack inequality.
\end{lemma}

\begin{remark}{\Cref{lem:li_yau_holds} is needed because \(W_d^s\) is not Lipschitz. By a suitable approximation argument, one should be able to extend the Li-Yau inequality systematically to general convex potentials and log-concave measures, without any boundedness or regularity requirements. The author has not been able to find a self-contained proof of this statement in the literature, and thus the proof proceeds through a slightly more elementary path.}
\end{remark}

\begin{proof}
    [Proof of \Cref{lem:li_yau_holds}]
    The Li-Yau inequality holds in \(\Barrel{\Omega}{V}\). For any \(f\in L^2(\Slice{\Omega}{V},\sqrt W_d)\), let \(f_F(x,w):=f(x,\sqrt{d}{2}-|w|)\). Now, by construction of \((\Slice{\Omega}{V},\sqrt d W_d)\), 
    \begin{equation}
    (P_t^{\Barrel{\Omega}{V}, 0} f_F) (x,w)
    =
    P_t^{\Slice{\Omega}{V},\sqrt d W_d}f(x,\sqrt{d}/{2}-|w|)        
    \end{equation}
    and thus one can extend the Li-Yau inequality from \(\Barrel{\Omega}{V}\) to \((\Slice{\Omega}{V},\sqrt d W_d)\).
\end{proof}

\begin{lemma}[Uniform ultracontractivity in \notoc{\(\Omega_d\)}]\label{lem:uniform_ultracont}
    Let \((\Omega,V)\) be a convex pair, with \(\lambda(\Omega, V)<2\) and let \((\Slice{\Omega}{ V}, \sqrt d W_d)\) be the slice domain associated to \(\Omega\), with the semigroup \(P^d_t:=P_t^{\Slice{\Omega}{V},\sqrt d W_d}\), and invariant measure \(\mu^s_d= Z_d^{-1} \sqrt{d}\exp(-\sqrt{d}W_d(x,r)) dx dr \). Then there exists a constant \(K = K_{(\Omega,V)}\) such that for all \(d\ge 0\):
    \begin{equation}
        \|P_{8}^d f\|_{L^\infty(S_d(\Omega,V))}\le K
        \|f\|_{L^2(S_d(\Omega,V),\mu_d^s(\Omega,\mu_d))}
    \end{equation}
\end{lemma}
\begin{proof}
    By ultracontractivity, it is enough to show that the conditions of Lemma~\ref{lem:ultracontractivity} hold. That exists a constant \(M_{\Omega,V}\) such that
    \begin{equation}
        \int_{\tilde S_d(\Omega, V)} e^{1/4(|x|^2+r^2)} d\mu^s_d \le M_{\Omega,V}
    \end{equation} follows from a direct computation. For large \(d\), using \Cref{lem:slice_estimates}, and that for \(r\in [0,\sqrt{d})\) one has \(\frac{r^2}{4}-\sqrt{d}r \le -\sqrt{d}{r}/2\)
    \begin{align}
        \int_{\tilde S_d(\Omega, V)} e^{\frac14(|x|^2+r^2))} d\mu^s_d(x,r)
        \le &
        2 Z_{V, \Omega}^{-1} |\Omega| e^{\operatorname{diam}(\Omega)^2} \int_0^{\sqrt{d}} \sqrt{d} e^{r^2/4} e^{-\sqrt{d}r} dr
        \\ \le &
        \tilde M_{\Omega, V}(\sqrt{d} \int_0^{\sqrt{d}}e^{-\sqrt{d}r/2} dx
        \\ \le &
        M_{\Omega, V}
    \end{align}

    To see that \(\|P_8^{(d)} \exp(x^2)\|\le M_{\Omega,V}\) we use again that for any function \(f \in L^2(\mu_d)\) the solution \(u((x,w),t):B_d(\Omega,V)\times \R_{\ge 0}\) to
    \begin{equation}
        \begin{cases}
            \partial_t u = \Delta u & \text{ in } B_d(\Omega,V) \\
            \partial_n u =0 & \text{ in } \partial B_d(\Omega,V) \\
            u((x,w),0) = f(x,\sqrt{d}/2-|w|)
        \end{cases}
    \end{equation}
    is by construction \(P_t^{d}(f)(x,\sqrt{d}/2-|w|)\) (as in the proof of \Cref{lem:li_yau_holds}). By the parabolic maximum principle, showing that \(\|P_8 \exp(x^2)\|\le M\) is equivalent to finding a barrier function \(u((x,w),t):\Omega_d\times \R_{\ge 0}\) such that
    \begin{equation}
    \label{eq:barrier_request}
    \begin{cases}
        \partial_t b \ge \Delta u& \text{ in } B_d(\Omega,V) \\
        \partial_n b \ge 0 & \text{ in } \partial B_d(\Omega,V) \\
        b((x,w),0) \ge \exp(1/4(\sqrt{d}/2-|w|)^2) & \text{ for all } (x,w) \text{ in } B_d(\Omega,V) \\
        b((x,w),8) \le M_{\Omega,V} &\text{ for all } (x,w) \text{ in } B_d(\Omega,V). \\
    \end{cases}
    \end{equation}

    The existence of such barrier function is the content of \Cref{lem:barrier_function}.
\end{proof}

As a corollary we obtain
\begin{corollary}
    The functions \(\phi_d^S\) are uniformly bounded in \(L^\infty(S_d(\Omega, V))\). Therefore the functions \(\phi_d\) are uniformly bounded in \(\Omega_d\).
\end{corollary}

\begin{proof}
    Consider only \(d\ge d_0\) such that \(\lambda_d \le 2.5\).
    The function \(\phi_d\) is an eigenfunction of \(L^{(d)}_{\Omega, V}\) of eigenvalue \(\lambda_d\le 2.5\),  satisfying \(\|\phi_d\|_{L^2(\mu^s_d(\Omega, V))}=1\). By the uniform ultracontractivity inequality
    \begin{equation}
        \|\phi_d\|_\infty = \|e^{8\lambda_d}P_t^{(d)}\phi_d\|_\infty \le e^{20}K_{(\Omega,V)}\|\phi_d\|_{L^2(\mu^s(\Omega, V))} = e^{20}K_{(\Omega,V)}
    \end{equation}
\end{proof}

The deferred ingredient in the proof of \Cref{lem:ultracontractivity} is as follows:
\begin{lemma}[Barrier function estimate]%
\label{lem:barrier_function}
    There exists a constant \(C_{\Omega,V}\) such that the function 
    \begin{equation}
            b(w,x,t) = C_{\Omega,V} + 20t +10d^{-1} |w|^2+\exp\left(\frac{d}{8}\right)\left(\frac{2}{2+4t}\right)^{d/2} \exp\left(-\frac{|w|^2}{2+4t}\right).
    \end{equation}
    satisfies~\eqref{eq:barrier_request} for some constant \(M_{\Omega, V}\).
\end{lemma}

\begin{proof}

 By enlarging \(M_{\Omega, V}\) larger a posteriori if necessary, one can assume that \(d\) is larger than \(d_0\) for some \(d_0\) to be fixed. The function \(b\) is a solution of the heat equation: The term \(20t +10d^{-1} |w|^2\) directly solves the heat equation, and the remainder is a multiple of the fundamental solution of the heat equation. 

 The function \(b\) does not depend on \(x\), and thus \(\partial_n b\) is zero in all of \(\partial B_{d}(\Omega,V)\cap (\partial\Omega \times \R^d)\). In the remaining of \(\partial B_{d}(\Omega,V)\), that is, in the set of points \begin{equation}\left\{(x,w), x\in \Omega, |w|= \frac{1}{2}\left(\sqrt{d}-\frac{V(x)}{\sqrt{d}}\right)\right\}\end{equation}
    the sign of \(\partial_n b(w,x,t)\) is that of \(w\cdot \grad_w b(w,x,t)\), which is equal to
    \begin{equation}
        2|w|^2 \left (10d^{-1} - \exp\left(\frac{d}{8}-\sqrt{d}\right)\left(\frac{2}{2+4t}\right)^{d/2} \exp\left(-\frac{|w|^2}{2+4t}\right) \right) =  2|w|^2 (*).
    \end{equation}
    The function 
    \begin{equation}
        d \left (\frac{1}{8} - \frac{1}{4}\frac{1}{2+4t} - \frac 12 \log\left(1+2t\right)\right)
    \end{equation}
    takes its maximum at \(t=0\) (which is equal to \(0\)). In particular, for \(|w|\ge \sqrt{d}/{2}\), if \(d\) is large enough,
    \begin{equation}
        (*) \ge   10d^{-1} - \exp(-\sqrt{d})\ge 0.
    \end{equation}
    This equation holds using the explicit constant \(10\), but any large constant instead of \(10\) would have sufficed for the proof.

    At \(t=0\), one needs to show that \(b(w,x,0)\ge \exp(1/4(\sqrt{d}/2-|w|)^2)\) for all \(x,w\) in the domain. This follows by showing that for all  \(|w| = \in [0, \sqrt{d}/2]\), the inequality
    \begin{equation}
        \label{eq:barrier_t0_first}
            \exp\left(\frac 1 4 \left|w-\frac{d}2\right|^2\right) \le \exp(d/8)\exp(-w^2/2)
    \end{equation}
    holds. The remaining case, when \(w\in[\sqrt{d}/2, \sqrt{d}/2+\frac{\|V\|_{\infty}}{2\sqrt{d}}]\), can be absorved into \(C_{\Omega,V}\).
    
    Showing ~\eqref{eq:barrier_t0_first} is equivalent to showing that for all \(|w| \in [0, \sqrt{d}/2]\),
    \begin{equation}
        \label{eq:barrier_t0_second}
        \frac 1 4 \left|w-\frac{\sqrt{d}}2\right|^2 
            \le
           d/8-w^2/2
    \end{equation}
    Let \(\tau = \frac 2 {\sqrt {d}} w\). Then~\eqref{eq:barrier_t0_second} follows from showing that for \(\tau \in [0,1]\),
    \begin{equation}
        \frac 1 4 \left|\tau-1\right|^2 
            \le
           \frac{1}{2}-\tau^2/2
    \end{equation}
    which holds as long as \(\tau \in [-1/3, 1]\supset [0,1]\).
    At \(t=8\) the inequality
    \begin{equation}
        \exp\left(\frac{d}{8}\right)\left(\frac{2}{2+4 t}\right)^{d/2} = \exp((1/8-\log\sqrt{17}))^d\le 1
    \end{equation}
    implies the uniform upper bounds.
\end{proof}

\subsubsection{\notoc{Dilating \(\psi_d\)}}

To show~\ref{item:limit} (a), we define another rescaled version of \(\phi\). Let \(\phi^{dil}_d(x,h):= \phi^S_d(x,2h/\sqrt{d})\). The domain of \(\phi_d^{dil}\) is 
\begin{equation}D_d:=\{(x,\sqrt{d} r), (x,r) \in S_d(\Omega, V)\} = \{(x,h), x\in \Omega, h\in [V(x), d]\}.\end{equation}

Let \(\mu^{dil}_d\) be the pushforward of \(\mu^{s}_d\) by the map \((x,r)\mapsto (x, \sqrt{d}r)\).  The measure \(\mu^{dil}_d\) has probability density
\begin{equation}
    \frac{1}{\tilde Z_d}\1_{D_d} \exp(-\sqrt{d}W_d(2x/\sqrt{d})) = \frac{1}{\tilde Z_d}\1_{D_d}  \exp\left(d \log\left(1-\frac{x}{d}\right)\right).
\end{equation}
The measures \(\mu^{dil}_d\) converge strongly to a measure \(\mu^{dil}_\infty\) with support in the set \(D_\infty:=\{(x,h), x\in \Omega, h>V(x)\}\) with density \(Z_{\Omega,V}^{-1}\exp(-x)\1_{D_\infty}\). Since \(d \log\left(1-\frac{x}{d}\right)-x = O(x^2/d)\) in \([-d/2,d/2]\), the Radon-Nikodym derivative \(\frac{d\mu^{dil}_d}{d\mu^{dil}_\infty}\) converges uniformly to \(1\) in the set \(\{(x,h): x\in \Omega, V(x)\le h \le d^{1/4} \}\).

The function \(\phi^{dil}_d(x,h)\) is, up to a sign, the unique minimizer of the problem
        \begin{equation}
        \label{eq:dilated_rayleigh}
            \phi^{dil}_d \in \argmin_{\int \psi d \mu^{dil}_d = 0} \frac{\int_{\Omega^{dil}_d} |\grad_x \psi|^2 + \frac d4 |\grad_h \psi|^2 d\mu^{dil}_d(x,h)}{\int_{\Omega^{dil}_d} |\psi| d\mu^{dil}_d(x,h)}.
        \end{equation}
The minimum value of the Rayleigh quotient problem is precisely \(\lambda_d\), and in particular, is uniformly bounded far from zero and infinity.

\begin{proposition}
    The functions \(\phi^{dil}_d(x,h)\) converge (up to a sign) as \(d\to \infty\) in \(L^2\) to a function constant in \(h\), and equal to the ground eigenfunction of \(L_{\Omega, V}\) in the \(x\) variable. 
\end{proposition}

\begin{proof}

Let \(w\) be a cutoff that is \(1\) for \(h\in [-1/4, 1/4]\) and \(0\) for \(|h|>\frac 12\). Then, as \(d\to \infty\) the functions \(\tilde \phi^{dil}_d(x,h)\) satisfy

\begin{equation}
\label{eq:rayleigh_weighted}
    \frac{\int_{\Omega^{\infty}} (|\grad_x  \phi^{dil}_d|^2 + d/4 |\grad_h  \phi^{dil}_d|^2) w(h/d^{1/4})d\mu^{\infty}(x,h)}{\int_{\Omega^{\infty}} | \phi^{dil}_d| d\mu^{\infty}(x,h)} \le \lambda_{\Omega_d}+o(1) =\lambda_{\Omega,V}+o(1)
\end{equation}

The functions \(\tilde \phi^{dil}_d(x,h)\) are uniformly bounded in \(W^{1,2}(w(h/d^{1/4})\mu_{\infty})\cap L^{\infty}(D_d)\). Therefore can extract a strongly convergent subsequence in \(L^2(\mu_{\infty})\) that is locally weakly convergent in \(W^{1,2}(\mu_{\infty})\). Let \(\hat \phi\) be the limit of this subsequence. Then

\begin{enumerate}
    \item \(\hat \phi\) has mean zero with respect to \(\mu_\infty^{dil}\).
    \item \(\int_{\Omega^\infty} |\hat \phi|^2 d\mu^{\infty} = 1\).
    \item By sending \(d\to \infty\) in~\eqref{eq:rayleigh_weighted},  \(\|\nabla \hat \phi\|_{L^2(\mu^{dil}_{\infty})} = 0\), and therefore \(\hat \phi\) is constant in \(h\). In other words, \(\hat \phi(x,h) = \hat \psi_1(x)\).
\end{enumerate}

By Fatou,
\begin{align}
    {\int_{\Omega^{\infty}} 
        |\grad_x \hat \phi|^2 d\mu^{dil}_{\infty}(x,h)}
    \le &
    \lim_{d\to \infty} {\int_{\Omega^{\infty}} 
        |\grad_x  \phi^{dil}_d|^2 + d |\grad_h  \phi^{dil}_d|^2 d\mu^{dil}_{\infty}(x,h)} 
    \\=& 
    \lambda_{\infty} \int_{\Omega^{\infty}} | \hat \phi|^2 d\mu^{dil}_{\infty}(x,h).
\end{align}
Since \(\hat \phi\) is constant \(h\) variable, integrating along it one obtains
\begin{equation}
    \begin{cases}
        \int_{\Omega} |\grad \hat \phi|^2 d\mu_{\Omega,V} \le \lambda_\infty
    \int_{\Omega} | \hat \phi|^2 \exp(-V) d\mu_{\Omega,V}\\
    \int_\Omega \phi \exp(-V) dx = 0.
    \end{cases}
\end{equation}
By uniqueness of the ground eigenfunction of \((\Omega,V)\), \(\hat \phi = \pm \phi_{\Omega,V}\). By changing the sign of each \(\phi_d\) as necessary, we see that \(\phi_d^{dil}(x,r)\to \phi_{\Omega,V}(x)\).
\end{proof}

From here one deduces \Cref{prop:equicontinuity}

\begin{proof}[Proof of \Cref{prop:equicontinuity}]
    The functions \(\phi^{S}_d\) are continuous in \(\overline{\Slice{\Omega}{V}}\), because \(\phi_d\) are continuous in \(\overline \Omega_d\). It suffices to find a family of functions \((f_d)_{d\ge 0}\), with \(f_d:\overline{\Slice{\Omega}{V}}\to \R\) that is equicontinuous, such that \(\|\phi^{S}_d-f_d\|_{L^{\infty}(S_d)}\) goes to zero.

    We will consider \(f_d:=\phi_{\Omega,V}^d(x, s):= e^{8\lambda_d}P_8^{\Slice{\Omega}{V}, W^s_d}  \phi_{\Omega,V}(x)\). The functions \((x,t)\mapsto \phi_{\Omega,V}(x)\) are equicontinuous in \(\Slice{\Omega}{V}\), and therefore by \Cref{thm:non-explosiveness}, so are the functions \(e^{8\lambda_d} P_8^{S_d, V}  \phi_{\Omega,V}(x)\) (note that \(\lambda_d\) are uniformly bounded).

    By ultracontractivity, 
    \begin{equation}
       \| e^{8\lambda_d} P_8^{\Slice{\Omega}{V}, W^s_d}  \phi_{\Omega,V} - \psi_d^S(x) \|_{L^\infty(\Slice{\Omega}{V})} \le C_{\Omega} e^{8\lambda_d} \|\phi_{\Omega,V}-\psi^{S}_d(x)\|_{L^2(\Slice{\Omega}{V}, W^s_d)} \to 0
    \end{equation}
\end{proof}

%% file: Parts/heat_flow.tex
Through this section we will consider the sequence of functions \(\psi^{com}_d(x,t)\) defined as \(\phi_{\Omega_d}(x,w) = \psi^{com}_d\left (x,1-\left (\frac{|w|}{\sqrt{d}/2}\right)^2\right)\) with domain
\begin{equation}
   H_d(B,V):= \{(x,t)\in \Omega \times \R^d, 1-(1-d^{-1}V(x))^2\le t\le 1\}
\end{equation}
Buy the assumption that \(V(x)\le 0\), the domains \(H_d(B,V)\) contain \(\Omega\times [0,1]\). By \Cref{prop:equicontinuity}, the functions \((\psi^{com}_d(x,y))_{d\ge 0}\) converge uniformly to \(\phi_{\Omega,V}(x)\) in the domains \((\{(x,t)\in H_d(B,V),t\le 0 \})_{d\ge 0}\).

By~\eqref{eq:psi_com}, To show~\ref{item:limit} of \Cref{thm:cylindrical_convergence} it suffices to show that the functions \(\psi^{(com)}_d\) are equicontinuous and converge to a solution of the heat equation in \(\Omega\times[0,1]\).

The functions \(\psi^{com}_d(x,t)\) are a (\(d-\)dependent) change of variable of \(\phi^S_d(x,s)\) in the second variable only. Therefore, by \Cref{prop:equicontinuity}, the functions \(\psi^{com}_d(x,t)\) are equicontinuous in space, and uniformly bounded, but a priori have no uniform continuity on the time variable. The continuity in the time variable will come from the explicit construction of the heat evolution of radially symmetric functions in \(\Barrel{\Omega}{V}\).

\begin{lemma}\label{lem:generic_eigenvalue}
    Let \(x_0 \in B_d(\Omega, V)^\circ\).  \(B^d_t(x_0)\) be a \((d+k+1)-\)dimensional Brownian motion in \(\Omega \times \R^{k+1}\) with reflective boundary conditions, and \(\tau\) a predictable stopping time that stops in \(B_d(\Omega,V)\) with \(\E[e^{2\lambda_d \tau}]<\infty\). Then 
    \begin{equation}
        \phi(x_0) = \E\left[e^{\lambda_d \tau} \phi(x_0+B_{\tau})\right].
    \end{equation}
\end{lemma}
\begin{proof}
    The random variable \(e^{\lambda_d(t\wedge \tau)}\phi(x_0 +B^d_{t\wedge \tau})\) is an \(L^{2}\)-martingale.
\end{proof}
After an Ito change of variables, one obtains the following result
\begin{lemma}
\label{lem:hitting_time_representation}

    For \(t\in (0,1]\), let \(H_s^{(t)}\) be a stochastic process satisfying 
    \begin{equation}
    \label{eq:SDE}
    \begin{cases}
        \dd H_s^{(t)} = -8(1+d^{-1}) \dd s + \frac 1 {d (1-H_s) } \dd B_s \\
        H_0^{(t)}=t\\
    \end{cases}
    \end{equation}
    and let \(\nu_d^{(t)}(s)\) be the law of the first time \(s^*\) at which \(H^{(t)}_{s^*}\) hits \(0\). Then for all \(d\) large enough
    \begin{equation}
        \label{eq:exponential_bounded}
        \E[\exp(2 \lambda_d s^*)] = \int_0^\infty e^{\lambda_d s} d\nu_d^{(t)} \le  1+ t.
    \end{equation}
    Let \(K^{\Omega}(x,x';t):\Omega\times \Omega \times [0, \infty)\to \R_{\ge 0}\) be the kernel of the heat flow in \(\Omega\) with Neumann boundary conditions. Then, for all \(d\) large enough
    
    \begin{equation}
    \label{eq:boundary_representation}
        \psi^{com}_d(x,t) = \int_{\Omega} \psi^{com}_d(x',1)  \int_0^\infty  K^{\Omega}(x,x';s) e^{\lambda_d s} d\nu_d^{(t)}(s) dx'.
    \end{equation}
\end{lemma}

\begin{proof}

    If \(X_t\) is a \((d+1)-\)dimensional Brownian motion, then \(1-\frac{\|X_t\|^2}{\sqrt{d}/2}\) satisfies~\eqref{eq:SDE}. This implies that once~\eqref{eq:exponential_bounded} is established,~\eqref{eq:boundary_representation} follows from  \Cref{lem:generic_eigenvalue}. In order to show that \(\E[\exp(2\lambda_d s^*)]\) is bounded, consider the process \(Y_s = (1+H_{s\wedge s^*}^{(t)})e^{2\lambda_d (s\wedge s^*)}\). We have
    \begin{equation}
        d Y_s = e^{2\lambda_d s} \left (  (2\lambda_d - 8(1+d^{-1})) \dd s + \frac 1 {d(1-H_s)}  \dd B_s\right)\1_{s<s^*}
    \end{equation}

    Note that \(\lambda_d \to \lambda_{\Omega, \infty}\le 2\), and in particular for \(d\) large enough, \(  (2\lambda_d - 8(1+d^{-1}))\le 0\). This shows that \(Y_s\) is a nonnegative supermartingale.  Moreover, \(\lim_{s\to \infty} Y_s =e^{\lambda_d s^*}\). In particular \(\E[e^{\lambda_d s^*}]\) is bounded by its value at zero, which is \(1+t\).
\end{proof}

By the same argument, for \(t>t'\), the measure \(\nu_{t\to t'}\) that characterizes the probability law of the hitting time of \(t'\) of \(H_s^{(t)}\) has the bound
\begin{equation}
    \mathbb E_{s\sim \nu_{t\to t'}}[\exp(2\lambda_ds)]\le 1+(t-t')
\end{equation}

From here, for example, one can extract uniform \(C^\alpha\) regularity in the spatial variable:
\
\begin{lemma}
    There are constants \(d_\Omega,M_\Omega>0\) such that for all \(d>d_\Omega\),
    for all \(0\le t \le t' \le 1\) with \(|t-t'|>\frac 1 {10}\), \(r\ge 0\) one has 
   \begin{equation}\label{eq:quantitative_spacetime}
       |\psi^{com}_d(x,t)-\psi^{com}_d(x,t')|\le \sup_{x'\in B_{C |t-t'|^{1/2}}(x)\cap \Omega} |\psi^{com}(x,t)-\psi^{com}(x',t)|+M_{\Omega} \left(\frac1{C^2}+|t-t'|\right)\|\psi_d^{com}\|_\infty
   \end{equation}
   In particular,  \(\|\psi^{com}_d(x,t)\|\)  are equicontinuous in \(\Omega \times [0,1]\).
\end{lemma}

\begin{proof}
     Let \(\nu_{t\to t;}\) be the law of the stopping of the first time \(H_s^{(t)}\)  hits \(t'\). We know that 
    \(\E_{s \sim \nu_{t'\to t}}[\exp(2\lambda_d s)] \le |t-t'|\), and thus \(\E_{s \sim \nu_{t'\to t}}[(\lambda_d s+1) \exp(\lambda_d s)] \le 1 +|t-t'|\). Rearranging,
    \begin{equation}
    \E_{s \sim \nu_{t'\to t}}[s \exp(\lambda_d s)] \le \lambda_d^{-1}|t-t'|.
    \end{equation}

    The kernel of the heat operator in  \(\Omega\) is the law of the reflected Brownian motion on a convex set. In particular, the boundary reflection further compresses \(K^\Omega(x,x',s)\) inwards (See \Cref{lem:heat_variance}):
    \begin{equation}
        \int K^{\Omega}(x,x',s) (x-x')^2 dx'\le
    \int K^{\R^{\operatorname{dim}(\Omega)}}(x,x',s) (x-x')^2 dx' = 
     \text{dim}(\Omega) \cdot s.
    \end{equation}

    This shows

    \begin{align}
        \int_{\Omega\setminus B_{C|t-t'|^{1/2}}(x)}   \int_0^\infty K^{\Omega} & (x,x';s) e^{\lambda_d s} d\nu_d^{(t)}(s) dx'\le
        \\\le &
        \frac{1}{C^2|t-t'|}
        \int_{\Omega}   \int_0^\infty  K^{\Omega}(x,x';s) (x-x')^2e^{\lambda_d s} d\nu_d^{(t)}(s) dx'
        \\ \le &
        \frac{\dim{\Omega}}{C^2|t-t'|}
        \int_{\Omega}   \int_0^\infty  K^{\Omega}(x,x';s) s e^{\lambda_d s} d\nu_d^{(t)}(s) dx'
        \\ = &
        \frac{\dim{\Omega}}{C^2|t-t'|}
          \int_0^\infty s e^{\lambda_d s} d\nu_d^{(t)}(s) dx'
        \\\le &
        \frac{\dim{\Omega}}{C^2\lambda_d^{-1}}.
    \end{align}

    For \(d_\Omega\) large enough so that \(\lambda_{\Omega,V}> \frac 12 \lim_{d\to \infty} \lambda_d\), one can choose \(C_\Omega = \frac{4 \dim{\Omega}}{C^2\lambda_{\Omega,V}}\). Combining this with
    \begin{equation}
        1\le \int_{\Omega}   \int_0^\infty  K^{\Omega}(x,x';s) e^{\lambda_d s} d\nu_d^{(t)}(s) dx = \int_0^\infty  e^{\lambda_d s} d\nu_d^{(t)}(s) dx \le 1 +|t-t|'
    \end{equation}
    gives the result.

    The functions \(\|\psi^{com}_d(x,t)\|\) are a change of variables in time of the functions \(\psi^{S}_d\), which are equicontinuous. Therefore, they are equicontinuous in space. By~\eqref{eq:quantitative_spacetime}, this equicontinuity is transferred to the time variable.
\end{proof}

The functions \(\psi_d\) are equicontinuous and satisfy Neumann boundary conditions on \(\partial \Omega \times [0,1]\). Then the last to show of \Cref{thm:cylindrical_convergence} (C.3) is.

\begin{proposition}
    Let \(\psi^{com}_\infty\) be a subsequential \(C^0\) limit of the \(\psi^{com}_d\). Then \(\psi_\infty\) solves the heat equation stated in \Cref{thm:cylindrical_convergence}\ref{item:limit}.
\end{proposition}

\begin{proof}
    By uniqueness of weak solutions to the heat equation (which for smooth \(\Omega\) follows by expanding a weak solution on an eigenbasis of \(\Delta\) in \(\Omega\)), it suffices to show that \(\psi^{com}_\infty\) satisfies the weak equation
    \begin{equation}
         8 \int_{\Omega} \psi^{com}_\infty \partial_t b
        + \int_{\Omega} \lambda_\infty \psi^{com}_\infty b dx + \int_\Omega \Delta_x b \psi^{com}_\infty  = \int_{\partial \Omega \times [0,1]} \grad_xb \cdot  \vec n ds
    \end{equation}
    for any \(C^2\) function \(b\) compactly supported on \(\R^k\times(0,1)\). But this follows from the  \(C^0\) convergence, since 
        \begin{equation}
        8\int_{\Omega} \psi^{com}_d \partial_t b
        +  \lambda_d \int_{\Omega} \psi_d^{com} b dx + \int_\Omega \Delta b \psi^{com}_d  = \int_{\partial \Omega \times [0,1]} \grad_xb \cdot  \vec n ds +(*)
    \end{equation}
where 
\begin{equation}
        (*) = \frac{1}{d} \int   \psi^{com}_d(x,t) (16 \partial_{tt} (b(x,t) (1-t)) - 6\partial_t b(x,t)) dx dt
    \end{equation}
\end{proof}

%% file: Parts/perturbation.tex
The goal of this section is to show \Cref{prop:perturbation}, which we rephrase below, recalling that \(R:=[-\pi/2,\pi/2] \times [-\pi/4,\pi/4]\).

\begin{proposition*}[\Cref{prop:perturbation} restated]
	Let \(q(y):[-\pi/4, \pi/4]\to \mathbb C\) be a smooth symmetric function with \(q'(\pm \pi/4) = 0\). There exists a smooth function \(\beta:R\to \mathbb R\) such that
	\begin{enumerate}
		\item \(\beta(x,y)\) is symmetric in \(y\),
		\item \(\beta(x,y)\) is antisymmetric in \(x\),
		\item \(\beta(\pi/2, y) - q(y)\) is constant,
	\end{enumerate}

	and a smooth convex potential \(V:R \to \mathbb R_{\le 0}\) such that as \(\epsilon \to 0\) the pair \((R,\epsilon V)\) has a spectral gap and
	\begin{equation}
		\begin{cases}
			\lambda_{R, \epsilon V} = 1 + O(\epsilon)\\
		\|\phi_{R, \epsilon V} - 2/\pi\sin(x)-\epsilon \beta \|_\infty = o(\epsilon)
		\end{cases}.
	\end{equation}
\end{proposition*}

\begin{proof}	

\textbf{1. Showing that one can differentiate the eigenfunctions of \(-\Delta + \epsilon \grad V\cdot \grad\) in \(\epsilon\).} 

For a smooth \(V\) in a rectangle, if \(\phi_0\) is the first eigenfunction of \(-\Delta+\epsilon \nabla V\cdot \nabla\), then \(\exp(-\epsilon/2 V)\phi_0\) is the first eigenfunction (with the same eigenvalue) of the operator
\begin{equation}Q_{\epsilon V}:=-\exp(+\epsilon/2 V)\nabla \exp(-\epsilon V) \nabla \exp(+\epsilon/2 V).\end{equation}

The operators \((I+Q_{\epsilon V})^{-1}\), are an analytic family of compact operators in \(L^2\), and for \(\epsilon = 0\), the operator \((I+Q_{0 V})^{-1}\) has no repeated eigenvalues. In particular, one can formally differentiate the eigenvalue equation \((I+Q_{\epsilon V})^{-1} \phi_{\epsilon} = \lambda \phi_{\epsilon}\) with respect to \(\epsilon\).

Since the eigenfunctions of \((I+Q_{\epsilon V})^{-1}\) are the same as those of \(Q_\epsilon\), the same holds for the conjugate operator \(\exp(-\epsilon/2 V)Q_\epsilon\exp(+\epsilon/2 V) = -\Delta+\nabla V\cdot\nabla\). In other words, one can differentiate the equation
\begin{equation}
	(-\Delta + \epsilon \nabla V \cdot \nabla  ) \phi_{R, \epsilon V} = 
	\lambda_{R, \epsilon V} \phi_{R, \epsilon V} 
\end{equation}
with respect to \(\epsilon.\)

Let \(\phi_{R,\epsilon V} = \phi_{R,0}+ \epsilon  \beta(x,y)+ O(\epsilon^2)\), and \(\lambda_{R,\epsilon V} = \lambda_{R,0}+ \epsilon  \mu_{R,V}+ O(\epsilon^2)\).  Then \(\psi_{R,V}\) solves the PDE
\begin{equation}
\label{eq:linearized}
	-\Delta \beta + \nabla V \cdot \nabla \phi_{R,0} =  \mu_{R,V}\phi_{R,0}+\lambda_{R,0}\beta.
\end{equation}

This equation determines both \(\psi_{R,V}\) and \(\mu_{R,V}\), since it implies that 
\begin{equation}
\mu_{R,V}\phi_{R,0} -  \nabla V \cdot \nabla \phi_{R,0} = 
	(-\Delta - \lambda_{R,0}) \beta \in \langle \phi_{R,0}\rangle ^\perp.
\end{equation}

Moreover, the residual terms \(r_\epsilon =  \epsilon^{-2}(\phi_{R,V}- \epsilon  \psi_{R,V})\) and \(\nu_\epsilon =  \epsilon^{-2}(\lambda_{R,V}- \epsilon  \mu_{R,V})\)  solve the elliptic equation
\begin{equation}
\label{eq:residue_eq}
	(-\Delta+\epsilon \nabla V\nabla)(\epsilon^{-2} \phi_{R,0}+\epsilon^{-1} \beta  + r_\epsilon)
	=
	(1+\epsilon \mu_{R,V}+\epsilon^2 \nu_{R,V}^{\epsilon})(\epsilon^{-2} \phi_{R,0}+\epsilon^{-1} \beta  + r_\epsilon)
\end{equation}
with \(\nu^\epsilon_{R,V} = O(1)\), \(\|r^\epsilon_{R,V}\|_2= O(1)\). Expanding~\eqref{eq:residue_eq}, and substituting in both~\eqref{eq:linearized} and the fact that \(-\Delta \phi_{R,0} = \lambda_{R,0}\) gives
\begin{equation}
\begin{aligned}
	(-\Delta+\epsilon V)(r_\epsilon) + \nabla V\nabla \beta = &
	(1+\epsilon \mu_{R,V}+\epsilon^2 \nu_{R,V}^{\epsilon})r_\epsilon
	\\&+
	(\mu_{R,V}+\epsilon \nu_{R,V}^{\epsilon}) 
	 \beta  
	 \\&+
	 \nu_{R,V}^{\epsilon}\phi_{R,0}
\end{aligned} 
\end{equation}

Therefore \(\|r_\epsilon\|_\infty \lesssim 1 + \|r_\epsilon\|_2 \lesssim 1\).  This reduces the problem to finding \(\beta, V\) with the required properties, that satisfy the linearized equation~\eqref{eq:linearized}.

\textbf{2. Finding a \(V\) that prescribes the right derivative.}

Let \(\beta\) be any smooth function that is symmetric in \(x\), antisymmetric in \(y\) with Neumann boundary data, equal to \(q(y)\) when \(x\) such that \((-\Delta+1)\beta_0(x,y)=0\) whenever \(x=\pi/2\). Then there exists a smooth \(V_0\) satisfying 
\begin{equation}
	(-\Delta - 1) \beta_0  = - \frac{2}{\pi}\cos(x)\partial_x V
\end{equation}

Let \(s^{(1)}(x), \mu^{(1)}\) solve the ODE
\begin{equation}
	x \cos x = (\partial_{xx}+1) s^{(1)}(x) + \mu^{(1)} \sin(x)
\end{equation}
in \([-\pi/2, \pi/1]\) with Neumann boundary conditions. The solution exists for some \(\mu^{(1)}\) by the Fredholm alternative and is smooth. Let \(\beta(x,y) := \beta_0(x,y)+ M s^{(1)}(x)\), and \(V(x):= V_0 + \frac 1 2 M (x^2+y^2)\), for \(M\) large enough so that \(V\) is convex. They satisfy the equation 
\begin{equation}
	(-\Delta - 1) \beta =  - M \mu^{(1)} \sin(x) - \frac{2}{\pi}\cos(x)\partial_x V,
\end{equation}
which finishes the proof.

\end{proof}

%% file: Parts/wings.tex
The goal of this section is to show
\begin{proposition*}
	Let \((R, V)\) be a convex pair, with \(R\) as in \eqref{eq:domain} such that \(\phi_{R, V}\) is antisymmetric in the first variable and \(V\) is bounded.  Let \(g(x,y)\) be a \(C^2([-1,1]\times[0,1])\) function such that:
	\begin{enumerate}
		\item For each \(x\), \(g(x,y)\) is convex in \(y\).
		\item \(g(x,y)\) is symmetric in the \(y\) variable.
		\item There is \(\delta>0\) such that \(g(x,y)=0\) for all \(x<\delta\).
	\end{enumerate}

	Let \(F(x_0,y_0): [0,1]\times[-1,1] \to [-1,1]\) be the value at \(x=0\) of the solution to the ODE
	\begin{equation}
		\begin{cases}
			f'(x) = \partial_y g(f(x))\\
			f(x_0)=y_0
		\end{cases}
	\end{equation}

	Then there exists a sequence of convex pairs \((R_m, V_{g,m})\), with \(R_m:=[-\pi/2 -m, \pi/2+m]\times [-1,1]\) such that, as \(m\to \infty\):
	\begin{enumerate}
		\item The pair \((R_m, V_{g,m})\) has a spectral gap and the ground eigenvalues \(\lambda_{R_m, V_{g,m}}\) converge to \(\lambda_{R,V}\).
		\item The functions \(\phi_{R_m, V_{g,m}}(x,y)\) are antisymmetric in \(y\).
		\item In \(R\), the functions \(\phi_{R_m, V_{g,m}}\) converge uniformly to \(\phi_{R,V}\)
		\item The functions \(T_m(x,y):=\phi_{R_m, V_{g,m}}(m x+\pi/2, \pi y)\) converge uniformly on \([0,1]\times[0,1]\)
		to \(T_{R,V,g}(x,y) :=\phi_{R,V}(\pi/2,F(x,y))\).
	\end{enumerate}
\end{proposition*}

\begin{lemma}%
\label{lem:pot_construction}
	There exists a sequence \(g_m:[-\delta,1]\times[-1,1]\to \mathbb R_+\) of bounded potentials such that 
	\begin{enumerate}
		\item The functions \(g_m\) converge in \(C^{2}([-\delta,1]\times[-1,1])\) to \(g\).
		\item The functions \(g_m(x,y)+m^2x + \frac{1}{2}m x_+^2\) are convex.
		\item For each \(g_m\) large enough there exists a family \((V_{l,m})_{l\ge 0}\) of symmetric convex potentials such that:
		\begin{enumerate}
			\item There exists a constant \(C_1\) such that \(g_{l,m}(x,y) = V(x,y)\) for \(|x|<\pi/2-C_1/l\).
		\item For \((x,y)\) in \([0,1]\times[-1,1]\), 
		\begin{equation} \nabla V_{l,m}(mx+\pi/2, y) = l \begin{pmatrix}
			m&0\\0&1
		\end{pmatrix}\vec v_m(x,y).\end{equation}
		where
		\begin{equation}
		\vec v_m(x,y) =(1+m^{-1} x+m^{-2}\partial_x g(x,y),\ \  \partial_y g(x,y)))^t
		\end{equation}
		\end{enumerate}
	\end{enumerate}
\end{lemma}

\begin{proof}
	The functions \(g_m\) will be of the form
	\begin{equation}
		g_m(x,y):=g(x,y)+C m^{-1} x_+^4(2+y^2)	
		\end{equation}
	for a large constant \(C\). The function \(x_+^4(2+y^2)\), for \(x\ge 0\), has Laplacian equal to \(2x^2(x^2+6(2+y^2))\) and Hessian determinant equal to \(8x^6(6-5y^2)\). They are

	\begin{itemize}
		\item Both are nonnegative in \([0,1]\times[-1,1]\). In particular, the Hessian of \(g_m+\frac{1}2 \sqrt{m}x^2 + mx\) is nonnegative definite in \([-\delta,\delta]\times[-1,1]\), where \(g\) is exactly zero.
		\item Both strictly positive in \([\delta,1]\times[-1,1]\), and thus there is a uniform \(\lambda_0\) such that the Hessian of \(x_+^4(2+y^2)\) is lower bounded by \(\lambda_0 I\) in \([\delta,1]\times[-1,1]\).
	\end{itemize} 

	For \(x\ge \delta\), the Hessian of \(g(x,y)+m^2x+ \frac{1}{2} m x^2\) is
	\begin{equation}
		\begin{pmatrix}
			g_{xx}& g_{xy}\\
			g_{yx}& g_{yy} 
		\end{pmatrix}
		+
		\begin{pmatrix}
			m&0\\
			0 & 0
		\end{pmatrix}
	\end{equation}
	Its trace is bounded from below by \(m-2\|g\|_{C^2}\), and its determinant is bounded from below by \(-\|g\|_{C^2}^2\). For \(m\gg \|g\|_{C^2}+1\) this shows that the lowest eigenvalue is at least \(-\frac{\|g\|_{C^2}^2}{m}\). In particular, there is a constant \(C\) such that for all \(m\), the lowest eigenvalue is at least \(-\frac{C}{m}\). Then the function \(g(x,y)+C m^{-1} x_+^4(2+y^2)+\frac 1 2 m x^2+ m^2 x\) is convex.

	Let \(h_{l,m}(x,y) = l(g_m(x,y) + \frac{1}{2}\sqrt{m} x^2+mx+\|V\|_{L^\infty(R)})\). The potential \(V_{m,l}\) will be defined, for \(m>\delta^{-1}\) as
	\begin{equation}
		\begin{cases}
			h_{l,m}((|x|-\pi/2)/m,y) \text{ if } |x|>\pi/2\\
			\max(V(x,y), h_{l,m}((|x|-\pi/2)/m,y)) \text{ otherwise.}
		\end{cases}
	\end{equation}

\end{proof}

In \([0,1]\times [-1,1]\), the functions \(\psi_{m,l}(x,y) = \phi_{R_m, V_{m,l}}(mx+\pi/2,y)\) satisfy the PDE
\begin{equation}
	l \vec v_m \cdot \nabla \psi_{lm} = (\partial_{yy} +m^{-2}\partial_{xx} +\lambda_{R_m, V_{m,l}}) = 0.
\end{equation}
Formally, sending \(l\) to infinity should give the PDE \(v_m \cdot \nabla \psi_{m}\), where \(\lim_{l\to \infty}\psi_{ml} = \psi_{m}\). That would imply that \(\psi_{ml}\) is constant on the flow lines by \(\vec v_m\). The map \((x,y) \mapsto (mx+\pi/2,y)\), maps \([0,1]\times [-1,1] \to [\pi/2,\pi/2+M]\times [-1,1]\). Therefore, the value \(\psi_{ml}\) should be determined by the location at which the flow lines by \(\vec v_m\) leave on the \emph{left}, at \(x=0\). This intuition is captured in the following lemmas:

\begin{lemma}%
\label{eq:fields_nice}
	Let \(m\) be large enough. Let \(v_m\) be the vector field
	\begin{equation}
	v_m(x,y):=(1+m^{-1}x+m^{-2}\partial_x g_m, \partial_y g_m) 	
	\end{equation} 
	with domain \(D:=[0,1]\times[-1,1]\). For all \(x_0\) in \(D\), the flow by \(v_m\) starting \(x_0\) exits \(D\) (for negative times) in the set \(\{0\}\times[-1,1]\). These vector fields converge  in \(C^1\) to the vector field
	\begin{equation}
		v(x,y):=(1, \partial_y g).
	\end{equation}
	The flow lines of \(v\) are graphs and are precisely the solutions to the ODE \(\frac{dy}{dx} = g(x,y(x))\).
\end{lemma}

\begin{proof}
	The vector field \(v_m(x,y)\) is symmetric with respect to flipping \(y\to -y\). If \(m\) is large enough, the first coordinate is strictly positive (and converges to 1), and for \(y\ge 0\), the second component is nonnegative (because \(g\) is symmetric and convex in \(y\)). In other words, the flow looks like
	\begin{center}
	\includegraphics[width=.12\textwidth]{Figures/flow.pdf}
	\end{center}
	Since \(g_m\) converges to \(g\) in the \(C^2\) topology, \(v_n\) converge to \(g\) in \(C^1\), for which the flow lines solve
	\begin{equation}
		\begin{cases}
			\dot x(t) = 1\\
		\dot y(t) = \partial_y g(x(t),y(t))
		\end{cases}
	\end{equation}
	which lets us rewrite them as graphs of the form \(y(x)\) solving \(\frac{dy}{dx} = g(x,y(x))\).
\end{proof}

\begin{lemma}%
\label{lem:l_to_infty}
	As \(l\to \infty\) the following hold:

	\begin{enumerate}
		\item \(\lim_{l\to \infty} \lambda_{R_m, V_{l,m}} = \lambda_{R,V}\)
		\item For \(l\) large enough (depending on \(m\)), the pair \((R_m, V_{l,m})\) has a spectral gap.
		\item The functions \(\phi_{R_m, V_{m,l}}\) are equicontinuous, and \(\lim_{l\to \infty}\|\phi_{R_m, V_{l,m}}-\phi_{R,V}\|_{L^\infty(R)}=0\).
		\item For large \(m\), the functions \(\phi_{R_m, V_{l,m}}\) converge as \(l\to \infty\) to a function \(\phi_m\), 
		uniquely determined in \(R\setminus R_m\) by \begin{equation}v_m \nabla(\phi_m(m(x+\pi/2),y))=0\end{equation}
	\end{enumerate}
\end{lemma}

\begin{proof}
	Fix \(m\) through the proof, and let \(l\to \infty\). Let \(\phi_{m,l}:= \phi_{R_m, V_{m,l}}\), and \(\lambda_{m,l}:= \lambda_{R_m, V_{m,l}}\).

	\begin{enumerate}
	\item The measures \(\exp(-V_{m,l})\) converge strongly to the measure \(\exp(-V)\1_R\).
	Therefore, by testing the Rayleigh quotient against \(\phi_{R,V}(\operatorname{sign}(x)\max(|x|,\pi/2),y)\), we see that
	\begin{equation}
		\limsup_{l\to \infty} \lambda_{R_m, V_{m, l}} \le \lambda_{R,M}.
	\end{equation}

	For all \(l\), the eigenfunctions \(\phi_{l,m}\) have the same domain, which is bounded, and have uniformly bounded eigenvalue. In particular, by \Cref{lem:ultracontractivity}, they are uniformly bounded in \(L^\infty(\mathbb R_m)\).

	With \(C_1\) as in Lemma \Cref{lem:pot_construction} 3.(a), the functions  \(\phi_{l,m}\) satisfy 
	\begin{equation}
		\limsup_{l\to \infty} \int_{-\pi/2+C_1/L}^{ \pi/2+C_1/l} \int_{-1}^1 |\nabla \phi_{m,l}|^2 \exp(-V(x)) dx dy \le \limsup_{l\to \infty} \lambda_{m,l} 
	\end{equation}
	and 
	\begin{equation}
		\lim_{l\to \infty} \int_{-\pi/2}^{ \pi/2} \int_{-1}^1 |\nabla \phi_{m,l}|^2 \exp(-V(x)) dx dy = 1.
	\end{equation}
	\item The \((R,V)\) has a spectral gap, any \(L^2\) weak subsequential limit of the \(\phi_{m,l}\) will saturate it. This shows that the functions \(\phi_{m,l}\) have to converge in \(L^2(R)\) to \(\phi_{R,V}\). This shows as well that \((R_m, V_{m,l})\) must have a spectral gap: If, for arbitrarily \(l\) the pair \((R_m, V_{m,l})\) didn't have a spectral gap, different choices of \(\phi_{R_m, V_{m,l}}\) would lead to different limits.
	\item Let \(f_m:=\phi_{R,V}(\operatorname{sign}(x)\max(|x|,\pi/2),y)\). The functions \(
		f_m-\phi_{m,l}\)
	are uniformly bounded and converge in \(L^2(R)\) to zero. They converge \(L^2(R_m, \exp(-V_m))\) converges to zero as well. Then, by ultracontractivity, the difference between \(\phi_{m,l}\) and \(e^{\lambda_{m,l}} P_1^{(R_m, V_{m,l})}f_m\) goes to zero in \(L^\infty\). But by non-explosiveness of flows (\Cref{thm:non-explosiveness}), the functions \(e^{\lambda_{m,l}} P_1^{(R_m, V_{m,l})}f_m\) are equicontinuous. This shows that the functions \(\phi_{m,l}\) are equicontinuous, and their convergence to \(\phi_{R, V}\) is not just in \(L^2\) but in \(L^\infty\).
	\item The functions \(\psi_{m,l}(x,y):=\phi_{m,l}(\pi/2+mx,y)\) in \([0,1]\times[-1,1]\) satisfy the PDE
	\begin{equation}
	 \vec v_m \cdot \nabla \psi_{m,l} = l^{-1}(\partial_{yy} +m^{-2}\partial_{xx} +\lambda_{R_m, V_{m,l}})\psi_{m,l}
	\end{equation}
	and by 3. converge uniformly in \(\{0\}\times[-1,1]\) to \(\phi_{R,V}(\pi/2,y)\). By integrating against a test function supported on the interior of \([0,1]\times[-1,1]\) and taking limits, any \(C^0\) subsequential limit of \(\psi_{m,l}\) is a weak solution to
	\begin{equation}
		\vec v_m \cdot \nabla \psi_{m,l} = 0
	\end{equation}
	that is continuous up to the boundary. By \Cref{eq:fields_nice}, this determines the value of all subsequential limits from the value in \(\{0\}\times[-1,1]\). Since the boundary data, given by \(\phi_{R,M}(\pi/2, \ast)\) is \(C^1\), so is the weak solution.
\end{enumerate}
\end{proof}

\begin{proof}
	[Proof of \Cref{prop:wings}]

	From \Cref{lem:l_to_infty}, one can extract a sequence of \(l_n\) growing fast enough so that:
	\begin{enumerate}
		\item The pair \((R_m, V_{m,l})\) has a spectral gap, and first eigenvalues converge to that of \((R, V)\).
		\item For \(m\) large enough the functions \(\phi_{R_m, V_m}\) are antisymmetric in \(y\): By symmetry of the domain and potential, they are either symmetric or antisymmetric, and the former cannot because they converge to an antisymmetric function.
		\item  Since for each \(m\) the functions \(\phi_{R_m, V_{m,l}}\) converge uniformly to \(V_{R,\phi}\), one can choose \(l_m\) to grow fast enough so that the same holds for \(\phi_{R_m, V_{m,l_m}}\).
		\item By continuity of ODE solutions with respect to the ODE coefficients, the functions \(S_m=\phi_m(m(x+\pi/2))\) converge to the unique solution in \([0,1]\times [-1,1]\) of \(\nabla v \times \nabla \psi = 0\) with left boundary data \(\phi_{R, V}(\pi/2, y)\). Again, by choosing \(l_m\) to grow fast enough, the same holds for \(\phi_{m, l_m}\).
		\end{enumerate}
\end{proof}

%% file: Parts/appendix.tex
The appendix contains standard computations on elliptic PDE and convex analysis.

\subsection{Step by step computations}%
\label{sub:step_by_step}

We show~\eqref{eq:psi_com_pde} step by step:

We start from \(\phi_{\Omega_d}(x,w) = \psi^{com}_d\left(x,1-\left (\frac{|w|}{\sqrt{d}/2}\right)^2\right)\), and the fact that for \(f:\R^k \to \R^{d+1}\), one has \(\Delta f(x,|w|) = \Delta_1 f(x,w) + \partial_{22} f+ \frac{d}{|w|} \partial_2 f\). In our situation, starting from
\begin{equation}
    - \Delta \phi_{\Omega_d}(x,w) = \lambda_{\Omega_d} \phi_{\Omega_d}(x,w)
\end{equation}
we obtain
\begin{equation}
\label{eq:before_expanding_polar}
-\left( \Delta_x+\Delta_\rho + \frac{d}{\rho} \partial_\rho\right) \left [ \psi^{com}_d\left (x,1-\left (\frac{\rho}{\sqrt{d}/2}\right)^2\right)\right]
=
     \lambda_{\Omega_d} \psi^{com}_d\left (x,1-\left (\frac{\rho}{\sqrt{d}/2}\right)^2\right).
\end{equation}
We compute the terms on the left-hand side: 
\begin{equation}
\label{eq:expanded_lapl_x}
    \Delta_x \left [ \psi^{com}_d\left (x,1-\left (\frac{\rho}{\sqrt{d}/2}\right)^2\right)\right] =
    (\Delta_1 \psi^{com}_d)\left (x,1-\left (\frac{\rho}{\sqrt{d}/2}\right)^2\right)
\end{equation}
\begin{equation}
    \label{eq:expanded_grad_rho}
    \partial_\rho \left [ \psi^{com}_d\left (x,1-\left (\frac{\rho}{\sqrt{d}/2}\right)^2\right)\right] 
    =
    -\frac{8\rho}{d}
    (\partial_2 \psi^{com}_d)\left (x,1-\left (\frac{\rho}{\sqrt{d}/2}\right)^2\right)
\end{equation}
\begin{align}
    \label{eq:expanded_lapl_rho}
    \Delta_\rho \left [ \psi^{com}_d\left (x,1-\left (\frac{\rho}{\sqrt{d}/2}\right)^2\right)\right] 
    = &
    \frac{64\rho^2}{d^2}
    (\partial_{22} \psi^{com}_d)\left (x,1-\left (\frac{\rho}{\sqrt{d}/2}\right)^2\right)
    \\&
    -\frac{8}{d}
    (\partial_2 \psi^{com}_d)\left (x,1-\left (\frac{\rho}{\sqrt{d}/2}\right)^2\right)
\end{align}

Therefore the left hand side of~\eqref{eq:before_expanding_polar} is equal to

\begin{equation}
    -\left(
    \underbrace
        {\Delta_1 \psi^{com}_d}
        _{\eqref{eq:expanded_lapl_x}}
    \underbrace
        {-8 \partial_2 \psi^{com}_d}
        _{\frac{d}{\rho}\eqref{eq:expanded_grad_rho}}
    +
    \underbrace
        {\frac{64\rho^2}{d^2}
        \partial_{22} \psi^{com}_d
        -\frac{8}{d}
        \partial_2 \psi^{com}_d)}
        _{\eqref{eq:expanded_lapl_rho}}
    \right)
    = \lambda_{\Omega_d}\psi^{com}_d
\end{equation}
Rearranging, and renaming variables to \((x,t)\) so that \(\rho= \frac{\sqrt{d}}{2}\sqrt{1-t}\), and \(\Delta_1 \to \Delta_x\), and \(\partial_d\to \partial_t\) one recovers
\begin{equation}
    8\partial_t \psi_d^{com}(x,t) = \Delta_x \psi_d^{com}(x,t)+\lambda_{\Omega_d}\psi_d^{com}(x,t) +
    \underbrace{
    \frac{1} d
    \left(
        16(1-t)
            \partial_{tt} \psi_d^{com}(x,t) 
            - 8 \partial_t  \psi_d^{com}(x,t) 
    \right)
        }
        _{\text{formally }\to 0}.
\end{equation}

We are interested in \(\rho \le \sqrt{2}{d}(1+o_d(1))\), which corresponds to \(t\) in \([-o(1), 1]\).

\subsection{Well-posedness and non-explosiveness of solutions to contractive flows}%
\label{sec:contractive_flows}

\begin{lemma}[Variance of the heat Kernel]\label{lem:heat_variance}
    Let \(\Omega\subseteq \mathbb R^d\) be a convex domain with smooth boundary. For \(x\in \Omega^\circ\), let \(K_\Omega(t,x_0,x)\) be the kernel of the heat evolution operator in \(\Omega\) with Neumann boundary data. Then 
    \begin{equation}
        \int_{\Omega} (x-x_0)^2 K_{\Omega}(t,x_0,x) dx \le 
        d\cdot t = \int_{\R^d} (x-x_0)^2 K_{\R^d}(t,x_0,x) dx 
    \end{equation}
\end{lemma}

\begin{proof}
    The function \(K_\Omega(t, x_0, \cdot)\) is the law, at time \(T\), of the reflected Brownian motion in \(\Omega\). This is the process satisfying the Skorokhod problem
    \begin{equation}
        \begin{cases}
            dW_t = dB_t - d\vec L\\
            W_0 = x_0
        \end{cases}
    \end{equation}
    where \(d \vec L\) is supported on the set \(\{W_t \in \partial_\Omega\}\), pointing inwards (i.e \(d\vec L \cdot (W_t-x_0)\le 0\)). Therefore 
    \begin{equation}
        d (W_t-x_0)^2 = d\cdot dt + 2 (W_t - x_0) dW_t +   2 (W_t - x_0) d\vec L 
        \le dt + 2 (W_t - x_0) dW_t 
    \end{equation}
    from which the result follows. 
\end{proof}

These are the same ingredients as in the proof of \Cref{thm:non-explosiveness}:

\begin{proof}
    [Proof of \Cref{thm:non-explosiveness}]
    Assume that \(\Omega\) is convex and \(V\) is in \(C^2(\overline \Omega)\), so that the reflected SDEs with drift \(-\nabla V\) is well posed (see e.g.~\cite{tanaka1979stochastic}). For each \(x\in \Omega\), let \(X_t^{(x)}\) be the solution to the Skorokhod problem
    \begin{equation}
        \begin{cases}
            d X_t^{(x)} = -\grad V (X_t) dt + dB_t - dL^{(x)}_t\\
            X_0 = x\\
            X_t^{x} \in \overline \Omega
        \end{cases}
    \end{equation}
    that is reflected in \(\partial \Omega\). We are using the coupling on all \(X^{(x)}_t\) arising from using the same Brownian motion \(dB_t\) for all \(x\). A direct computation shows
    \begin{equation}
        d |X^{(x)}_t-X^{(y)}_t|^2 = - 2 (X^{(x)}_t-X^{(y)}_t)\cdot (\grad V(X^{(x)}_t)-\grad V(X^{(y)}_t)) - 2 (X^{(x)}_t-X^{(y)}_t) (dL^{(x)}_t - dL^{(y)}_t).
    \end{equation}
    Since \(V\) is convex, the first term is negative. Similarly, since \(\Omega\) is convex, for any \(x\in \Omega\) and any \(y\in \partial \Omega\), \((x-y)\cdot \vec n_y\ge 0\), where \(\vec n_y\) is the inwards-pointing normal. In particular, since \(dL^{(x)}_t\) is supported on the times when \(X^{(y)}_t\) is in \(\partial \Omega\), proportional to \(\vec n_{X^{(y)}_t}\) and pointing inwards, the second summand is also nonnegative. 

    In particular, \(d |X^{(x)}_t-X^{(y)}_t|^2 \le 0\). This shows that if \(\nu_{x,t}\) is the law of \(X^{(x)}_t\), then one has the \(\infty\)-Wasserstein bound
    \begin{equation}
        W_\infty(\nu_{x,t}, \nu_{y,t}) \le |x-y|.
    \end{equation}

    But by Feynman-Kac,
    \begin{equation}
        P_t^{\Omega,V}(f) (x)= \int f(x') d\nu_{x,t}(x')
    \end{equation}
    from which the result follows by Wasserstein duality.
\end{proof}

The regularity of \(V\) is not necessary, as seen by the following approximation argument:
\begin{lemma}\label{lem:fokker_plank_stable}
    Let \((V,\Omega)\) be a convex pair, such that \(\|\nabla V\|_{L^\infty(\Omega)} <\infty\).
    Let \(((V_n,\Omega_n))_{n\ge 0}\) be a sequence of convex pairs such that the sets \(\Omega_n\) are increasing and converge (in the Hausdorff sense) to \(\Omega\), and \(\|V_n-V\|_{L^\infty(\Omega_n \cap \Omega)} \to 0\). Then
    \begin{enumerate}
        \item \(Z_n:=\int_{\Omega n} \exp(-V_n) dx\) converges to \(Z:=\int_{\Omega n} \exp(-V_n) dx\).
        \item For every eigenvalue \(\lambda_{\Omega,V}^{(i)}\) (not just the first one), \(|\lambda_{\Omega_n,V_n}^{(i)}-\lambda_{\Omega,V}^{(i)}| \to 0\).
        \item If \((\Omega,V)\) has a spectral gap, for all \(n\) large enough, the pair \((\Omega_n, V_n)\) has a spectral gap as well.
        \item \(\|\phi_{\Omega_n, V_n}-\phi_{\Omega,v}\|_{L^\infty(\Omega_n \cap \Omega)} \to 0\).

        \item For any function \(f\in L^2(\Omega)\) and any \(t>0\), the functions \(P_t^{\Omega_i, V_i} f\) converge uniformly to \(P_t^{\Omega,V}f\).
    \end{enumerate}
\end{lemma}
\begin{proof} 
    Let \(\phi^{(i)}:= \phi^{(i)}_{\Omega,V}\), and \(\phi^{(i)}_n:= \phi^{(i)}_{\Omega_n,V_n}\) be  the eigenfunctions of \((\Omega, V)\) and \((\Omega_n, V_n)\). These functions are Holder-continuous for each \(n\), but a priori not uniformly.
    \begin{enumerate}
        \item Follows from dominated convergence. By redefining the potentials to be \(V_n-\log Z_n\), one may assume that \(Z_n=1\) for the proof of parts 2-5. 
        \item Let \(\phi^{(i)}\) be the \(i-\)th eigenfunction of the limiting pair \((\Omega, V)\). Then (by e.g. dominated convergence),
    \begin{equation}
        \lim_{n\to \infty} \int_{\Omega_n} \phi^{(i)}(x)\phi^{(j)}(x) \exp(-V_n(x))  dx\to \delta_{ij}
    \end{equation}
    and 
    \begin{equation}
        \lim_{n\to \infty} \int_{\Omega_n} \grad \phi^{(i)}\cdot\grad \phi^{(j)} \exp(-V_n(x))  dx \to \lambda^{(i)}\delta_{ij}.
    \end{equation}
    By the Courant-Fischer variational formula for eigenvalues this shows that,
    \begin{equation}
        \limsup_{n\to \infty} \lambda^{(i)}_n \le \lambda^{(i)}.
    \end{equation}

    By the Li-Yau inequality, each fixed \(i\), the functions \(\phi^{(i)}_n\) are then uniformly bounded 
    (see \Cref{lem:ultracontractivity}) and are uniformly bounded in \(H^1(\Omega_i)\). 
    In particular, the functions \(\phi^{(i)}_n\) (extended to be e.g. zero in \(\Omega_n^c\)) have strong subsequential limits in \(L^2(\Omega)\), and weak subsequential limits in \(H^1_{loc}(\Omega)\). By passing to a subsequence, one may assume that \(\phi^{(i)}_n\) converges to \(\phi^{(i)}_\infty\) in \(L^2(\Omega)\). Then

    \begin{equation}
        \int_{\Omega} \phi^{(i)}_\infty (x) \phi^{(j)}_\infty(x) dx = \delta_{ij}.
    \end{equation}

    Let \(\delta_n\to 0\) such that \((1+\delta_n)\Omega_n \supseteq \Omega_0\). Fix some value of \(n\ge 0\). The functions \(\tilde \phi_n^{(i)}:=\phi^{(i)}_n((1+\delta_n)^{-1}x)\) converge to \(\phi^{(i)}_\infty\). Let \(\psi_n:= \sum_{j=1}^i \alpha_i \tilde \phi_n^{(i)}\), for some coefficients \(\alpha_1, \dots \alpha_i\) such that \( \sum_{j=1}^i \alpha_i^2=1\), and \(\psi_\infty:= \lim_{n\to \infty} \psi_i\). By Fatou, 

    \begin{equation}
        \int_{\Omega} |\nabla \psi_\infty (x)|^2 dx \le \liminf_{n\to \infty} \lambda^{(i)}_n.
    \end{equation}
    By the variational inequality for eigenvalues, this shows that  \(\lambda^{(i)}_{\Omega,V}\le \lim_{n\to \infty} \lambda^{(i)}_{\Omega_n,V_n}\)
    \item Is a direct consequence of 2.

    \item Item 2. shows that \(\phi^{(i)}_\infty\) (which are the \(L^2\) limits of  \(\phi^{(i)}_n\) as \(n\to \infty\)) are a family of eigenfunctions of \((\Omega,V)\), with the right multiplicity. It suffices to show that this limits are in \(L^\infty\) as well.
    
    By elliptic regularity for bounded-measurable coefficients on Lipschitz domains, the functions \(\phi^{(i)}(\Omega, V)\) are \(C^\alpha\) Holder-regular, with constant \(C_0 \lambda^{(i)}(\Omega, V)\). The value of \(\alpha, C_0\) depend only on the dimension,  the Lipschitz connstant of \(\Omega_i\), and the norm \(\|\nabla V_i\|_\infty\). This shows that for each \(i\), the eigenfunctions \(\phi^{(i)}_{n}\) are equicontinuous, and thus converge to \(\phi_\infty^{(i)}\) uniformly and not just in \(L^2\).

    \item

    The \(L^2\) convergence of all Laplace eigenfunctions implies that, for any \(f\in L^2\) and any \(t>0\), the functions \(P_t^{\Omega_n, V_n} f\) converge in \(L^\infty_{{[t, \infty)}}L^2_\Omega\) to \(P_t^{\Omega, V} f\) as \(n\to \infty\). It remains to show that the convergence is in \(L^\infty_\Omega\).

    Let \(\epsilon>0\), \(t>0\). There is a \(j_0\) and \(n_0\) such that for all \(n>n_0\) (including \(n=\infty\), with \(\Omega_\infty = \Omega\), \(V_\infty = V\)), \(j\ge j_0\)  and any function \(f \in L^2(\Omega_n,V_n)\) such that if \(\langle \phi^{(j)}_n, f\rangle_{L^2(\Omega_n, \exp(-V_n))} =0\), one has the estimate
    \begin{equation}
        \|P_{t/2}^{\Omega_n,V_n}f\|_{L^2(\Omega_n, \exp(-V_n))} \le \epsilon \|f\|_{L^2(\Omega_n, \exp(-V_n))}.
    \end{equation} 
    By Li-Yau, this shows that if \(j_0,n_0\) are large enough,  one can ensure
    \begin{equation}
        \|P_{t}^{\Omega_n,V_n}f\|_\infty \le  \epsilon \|f\|_{L^2(\Omega_n, \exp(-V_n))}. 
    \end{equation}

    Now by triangle inequality
    \begin{equation}
        \|P_t^{\Omega_n, V_n} f - P_t^{\Omega,V} f \|_{L^{\infty}(\Omega_n)} = \sum_{j\le j_0}
        \underbrace{
        \|\phi^{(j)}_ne^{-\lambda^{(j)}_n t}\langle \phi^{(j)}_n, f\rangle - \phi^{(j)}e^{-\lambda^{(j)} t}\langle \phi^{(j)}, f\rangle  \|_{L^\infty}}_{\text{ Vanishes by 4. as } n \to \infty}
        + 2 \epsilon \|f\|_{L^2(\Omega)}
    \end{equation}
    \end{enumerate}
\end{proof}

While the approximating argument above still requires \(\nabla V\) and the \(\nabla V_n\) to be uniformly bounded, the Lemma above should hold with just pointwise convergence to a (potentially unbounded) \(V\).

\subsection{Poincaré constants}
\begin{lemma}
\label{lem:poincare_ball}
    The Poincaré constant of \(B_1(\R^{d+1})\) is at most \(d^{-1}\). Equivalently, \(\lambda_{B_1(\R^{d+1})}\ge d\).
\end{lemma}

\begin{proof}
    
    The lowest eigenvalue of \(\mathbb S^{d}\), the \(d\)-dimensional sphere of radius \(1\) is the eigenvalue of the lowest spherical harmonic, which is \(d\).

    Let \(\psi_d\) be the ground eigenfunction of \(B_1(\R^d)\). There are two options. If for every \(\rho \le 1\) the function \(\int_{\mathbb S^{d-1}} \psi_d(\rho \theta) d\theta = 0\), by the sharp Poincaré inequality on spheres
    \begin{align}
    \int_{B_1(\R^d)} |\psi_d(\rho \theta)|^2 d\theta  = &
    \int_0^1 \rho^{d}
        \int_{\mathbb S^{d-1}} |\psi_d(\rho \theta)|^2 d\theta d\rho 
    \\\le &
    \int_0^1 \rho^{d}
        (d^{-1}\rho^2)
        \int_{\mathbb S^{d-1}} |\grad_{\rho \mathbb S^{d-1} }\psi_d(\rho \theta)|^2 d\theta d\rho
        \\ \le &
        d^{-1}
  \int_0^1 \rho^{d-1}
        \int_{\mathbb S^{d-1}} |\grad \psi_d(\rho \theta)|^2 d\theta d\rho
        \\=&
    d^{-1}
    \int_{B_1(\R^d)} |\grad \psi_d(\rho \theta)|^2 d\theta  
    \end{align}
    and thus the ground eigenvalue is at least \(d\). If, on the other hand, there exists a \(\rho\) such that \(\int_{\mathbb S^{d-1}} \psi_d(\rho \theta) d\theta \neq 0\), the function \(\eta(\rho):= \fint_{\mathbb S^{d-1}} \psi_d(\rho \theta) d\theta \). Then the function \(\eta(|w|)\) is a non-zero Neumann eigenfunction of \(\Delta\) with the same eigenvalue as \(\psi_d\) in \(B_1(\R^{d+1})\). In this case, showing that the eigenvalue is at least \(d\) corresponds to showing the equality in
    \begin{equation}
            \int_{B_1(\R^d)} |\psi_d(\rho \theta)|^2 d\theta  =
            |\mathbb S^{d}|\int_0^1 \rho^{d}\eta(\rho) d\rho 
            \le 
            d^{-1}
            |\int_0^1 \rho^{d}|\eta'(\rho)|^2 d\rho
            = 
            \int_{B_1(\R^d)} |\psi_d(\rho \theta)|^2 d\theta
    \end{equation}

    The weight \(\rho^{d}\) is \(d-\)strongly log-concave in \([0,1]\) (in the sense that \(-\log(x^d)''\ge d\)), and therefore, by (e.g) the Bakry-Emery criterion, the Poincaré constant is at most \(d^{-1}\).
\end{proof}

\subsection{Ultracontractivity}

\begin{lemma}[Ultracontractivity (Following e.g. \cite{rockner2003supercontractivity})]
\label{lem:ultracontractivity_proof}
    Let \((\Omega,V)\) be a convex pair with measure \(\mu:=Z_{\Omega,V} \exp(-V) \1_\Omega\) and heat semigroup \((P_t)_{t\ge 0}\) for which the Li-Yau inequality holds. Assume that there exists a \(t\ge 0\) and \(C>0\) such that
    \begin{equation}
        \|P_{t/2}(\exp(2|x|^2/t))\|_{\infty} \le M \text{, and }
        \int_{\Omega} \exp(-2|x|^2/t) d\mu \ge M^{-1}
    \end{equation}
    then \(\|P_t f\|_{L^\infty(\Omega)} \le M^2 \|f\|_{L^2(\mu))}\).
    In particular, for any \(t>0\),
    \(\|P_t f\|_\infty \le C_{t, \operatorname{diam}(\Omega)} \|f\|_{\infty}\).
\end{lemma}

\begin{proof}
    The Langevin flow \(P_t\) satisfies a maximum principle, and so it suffices to show that for any \(f\) such that \(\|f\|_{L^2(\mu)}=1\), one has
    \begin{equation}
        |P_{t/2} f|(x) \le C \exp(|x|^2/t).
    \end{equation}

    By the Li-Yau inequality,
    \begin{equation}
        |P_{t/2}(f)|^2(x) \exp\left(-t^{-1}|x-y|^2\right) \le |P_{t/2}(f^2)|(y).
    \end{equation}
    The flow \(P_t\) preserves integrals against \(d\mu\), and therefore
    \begin{equation}
        |P_{t/2} f|(x) \int \exp\left(-t^{-1}|x-y|^2\right) d\mu(y) \le \int P_{t/2}f^2(y)d\mu(y) = 1.
     \end{equation} 
     But
     \begin{equation}
         \exp\left(-t^{-1}|x-y|^2\right) \le \exp(-2 t^{-1} |x|^2)\exp(-2t^{-1} |y|^2)
     \end{equation}
     and thus
     \begin{equation}
         |P_{t/2} f|(x) \le C \exp(2 t^{-1} |x|^2)
     \end{equation}

\end{proof}

\subsection{Maximum principle for the Neumann heat equation}

\begin{lemma}[Maximum principle for the Neumann heat equation]
\label{eq:neumann_max_smooth}
    Let \(\Omega\) be a bounded convex domain, and \(b\) be a \(C^1\) subsolution of the heat equation on \(\overline \Omega \times (0,1]\), that is
    \begin{equation}
        \begin{cases}
            \partial_t b \le \Delta u & \text{ in } \Omega ^\circ \times (0,1)\\
            \partial_n b \le 0 & \text{ in } (0,1) \times \partial \Omega.
        \end{cases}
    \end{equation}
    with \(u(x,0)\le 0\). Then \(u(x,t)\le 0\).
\end{lemma}

\begin{proof}
    Assume that \(0 \in \Omega^{\circ}\), then let \(b^\epsilon(x,t):=b(x,t)+2 \epsilon t + \frac 1 {2d} \epsilon |x|^2\). Then \(b_\epsilon\) is also a sub-solution, satisfying
    \begin{equation}
        \begin{cases}
            \partial_t b^\epsilon \le \Delta b^\epsilon -\epsilon & \text{ in } \Omega ^\circ \times (0,1]\\
            \partial_n b^\epsilon \le -C_\Omega \epsilon & \text{ in } (0,1] \times \partial \Omega.
        \end{cases}
    \end{equation}
    In particular, by the second equation, the function \(u^\epsilon\) cannot take its maximum in \(\partial \Omega \times [0,1]\). Similarly, it cannot be at \(t=1\), or on the interior. (Note that convexity allows us to bypass Hopf Lemma).
\end{proof}

We say that \(B\) is a super-solution if \(-B\) is a solution. Then

\begin{corollary}[Maximum principle with barriers]
    Let \(\Omega\subset \R^d\) be a convex, bounded domain. If \(u(x,t)\) is a sub-solution to the Neumann heat equation in \(\Omega\), \(B(x,t)\) is a super-solution, and \(b(x,t)\) is a sub-solution to the Neumann heat equation in \(\Omega\), such such that \(b(x,0)\le u(x,0)\le B(x,0)\) for all \(x\in \Omega\), then for all \(t\ge 0\)
    \begin{equation}b(x,t)\le u(x,t)\le B(x,t)\end{equation}
\end{corollary}
\begin{proof}
    By \Cref{lem:fokker_plank_stable}, one may assume that \(\Omega\) is smooth. In this case, the solution to the heat equation is \(C^1\) up to the boundary, and one can apply \Cref{eq:neumann_max_smooth} to \(b(x,t)-u(x,t)\) and to \(u(x,t)-B(x,t)\).
\end{proof}